\documentclass[a4paper]{scrartcl}

\usepackage{a4wide}
\usepackage[utf8]{inputenc}
\usepackage{graphicx}
\usepackage{hyperref}
\usepackage{float}
\usepackage{amssymb}
\usepackage{amsmath}
\usepackage{amsthm}
\usepackage{mathtools}
\usepackage{authblk}

\usepackage[ruled,vlined,linesnumbered]{algorithm2e}
\usepackage{extarrows,todonotes}
\usepackage{enumitem}

\usepackage{subfigure}

\newcommand{\OWA}{\textup{OWA}}
\newcommand{\VAR}{\textup{VaR}}

\newcommand{\X}{{\mathcal{X}}}
\newcommand{\cU}{{\mathcal{U}}}
\newcommand{\E}{\mathbb{E}}
\newcommand{\Sp}{\rho}
\newcommand{\V}{\mathbb{V}}
\newcommand{\R}{\mathbb{R}}
\newcommand{\N}{\mathbb{N}}

\newtheorem{theorem}{Theorem}

\newtheorem{definition}{Definition}

\newtheorem{lemma}{Lemma}

\newtheorem{property}{Property}
\newtheorem{cor}{Corollary}
\newtheorem{remark}{Remark}

\begin{document}

\title{An extension of Ordered Weighted Averaging  over intervals with application to optimization under risk}

 \author[1]{Werner Baak}
\author[1]{Marc Goerigk\footnote{Corresponding author.}}

\affil[1]{Business Decisions and Data Science, University of Passau, Germany,

 \texttt{\{werner.baak,marc.goerigk\}@uni-passau.de}}

\author[2]{Adam Kasperski}

\affil[2]{Department of Operations Research
and Business Intelligence, Wroc{\l}aw University of Science and Technology, Poland,

\texttt{adam.kasperski@pwr.edu.pl}}

\author[3]{Pawe{\l} Zieli\'nski}

\affil[3]{Department of Fundamentals of Computer Science, Wroc{\l}aw University of Science and Technology, Poland,

\texttt{pawel.zielinski@pwr.edu.pl}}

\date{}

\maketitle

\begin{abstract}
The Ordered Weighted Averaging (OWA) operator is a traditional and commonly used criterion for aggregating discrete values of uncertain quantities. In this paper, it is shown that the discrete OWA naturally extends to the continuous case by using the concept of a distortion risk measure. It is shown how to apply the distortion risk measure to optimization problems with a linear objective function, whose coefficients are random variables with continuous  distribution functions supported on intervals. The case where these coefficients are independent, uniformly distributed random variables is explored in more detail.  
The computational complexity of the resulting optimization problem is analyzed, and solution methods with approximation guarantees are proposed.
These methods are also verified through computational experiments.
\end{abstract}


\section{Introduction}

Decision-making under uncertainty is one of the fundamental challenges in operations research and management science. The Ordered Weighted Averaging (OWA) operator, introduced by Yager~\cite{YA88}, has emerged as a powerful tool to aggregate multiple scenarios in decision-making processes. OWA operators provide a flexible framework that can model a wide range of decision-maker attitudes, from optimistic to pessimistic, by adjusting the weighting vector~\cite{BGH24, XU05}.
The theory and applications of OWA operators have been extensively developed since their introduction (see, e.g.,~\cite{YKB11}).

The traditional discrete OWA operator proposed in~\cite{YA88} aggregates a finite collection of values, typically corresponding to a finite number of objective functions in decision-making models. An important special case includes decision making under risk, where the objective functions are induced by scenarios (possible states of the world). Using non-increasing weight vectors, OWA can capture a worst-case or robust optimization perspective \cite{BN09, KY97, GH24}. This risk-averse framework also aligns well with criteria such as Value-at-Risk (VaR) or Conditional Value-at-Risk (CVaR) \cite{bertsimas2009constructing, rockafellar2000optimization, P00}, which focus on limiting exposure to extreme losses beyond a specified confidence level. Optimization problems with the discrete OWA operator have been widely discussed in the literature. In~\cite{CG15, GS12, OO12, OS03} mixed integer programming formulations for the general problem and its special cases were proposed. The computational complexity of minimizing OWA in combinatorial problems was investigated in~\cite{KZ15}. Except for very specific weights, the problem turned out to be NP-hard even for simple selection problems (see~\cite{KZ15}). Finally, some approximation algorithms for minimizing OWA have been proposed in~\cite{BKKZ25,CGK20}.

Discrete uncertainty is one several methods of  uncertainty representations. The second, commonly used method is the interval uncertainty (see, e.g.,~\cite{KY97}),  in which closed intervals are assigned to uncertain parameters. Unfortunately, the discrete OWA operator is not suitable to handle such uncertainty, because the number of scenarios (realizations of the parameters) can be  infinite.
In the literature, some attempts have been made to extend OWA to the interval case. One, called C-OWA, has been proposed in~\cite{YA04} and further extended in~\cite{YX06, ZC10}. The idea is to aggregate the values in a closed interval  using a weight function defined on the interval $[0,1]$. In this paper, we will point out some limitations of the C-OWA operators when they are applied to optimization problems. 
It turns out that discrete OWA naturally extends to the continuous case by applying the concept of a distortion risk measure~\cite{D90, SBR10}. The risk measure is a function that assigns a real value to a random variable, typically representing a random loss (see, e.g.,~\cite{KAM11}). In the distortion risk measure, we use a weight function defined in $[0,1]$ to give more importance to the specified loss values (for example, to larger values). Therefore, the idea is similar as in the discrete OWA operator, where we assign weights to specified positions of ordered sequence of discrete values. Furthermore, if we apply the distortion risk measure to a single uniformly distributed random variable, then we get the C-OWA operator defined in~\cite{YA04}. 

In this paper, we wish to investigate the class of optimization problems with a linear objective function. We assume that the objective function coefficients are random variables with interval supports. In consequence, the objective function value is also a random variable with interval support.  We will consider a general case and, more deeply, a model, in which the objective function coefficients are independent and uniformly distributed random variables. In robust optimization, this corresponds to the interval uncertainty representation, where the uncertainty set is a Cartesian product of closed intervals. In previous approaches discussed in the literature, the min-max (regret)  criterion~\cite{KY97}, or its budgeted version~\cite{BS04} was used to compute a solution.
 In this paper, we present an alternative and more sophisticated model in which a distortion risk measure is applied. Our approach has two advantages. Firstly, it captures the fact that  the values of the objective function are not equally probable, even if all coefficients have uniform probability distributions. Notice that in robust optimization, we usually ignore the probabilistic information by assuming that the probability distribution for the objective function values is unknown. This leads to application of the simple min-max criteria that take into account only the worst scenarios. Secondly, we can model decision-maker risk aversion by specifying a weight function. We will assume that this weight function can be generated from a basic unit monotone (BUM) function~\cite{YA04}. If the BUM function is concave, then the distortion risk measure is coherent~\cite{SBR10}, which implies some favorable computational properties of the resulting optimization model that we will exploit.

This paper is organized as follows. In Section~\ref{sec1}, we recall the definition of the traditional discrete OWA operator and show its interpretation in a stochastic setting. We then recall the definition of distortion risk measure and its special case, called spectral risk measure. We show that distortion risk measure is a natural extension of discrete OWA to the continuous case. We prove several properties of the distortion risk measure that will be used later in the optimization context. In Section~\ref{sec:2} we consider an optimization problem with uncertain objective function coefficients. We apply distortion risk measure to compute a solution to this problem and propose some general solution algorithms. We then investigate a special case in which the objective function coefficients are independent uniformly distributed random variables. We prove that evaluating a given solution in such a problem is $\#$P-hard. This fact motivates us to construct some approximation algorithms for the problem. For the concave BUM function, which induces a non-increasing weight function, we are able to show some approximation guarantees for these algorithms. Section~\ref{sec4} contains the results of some computational tests.

\section{Discrete OWA operator and distortion risk measure}
\label{sec1}

In this section we first recall the definition of the traditional  discrete OWA operator proposed in~\cite{YA88} and provide its interpretation in the stochastic setting. Next, we  show that the well-known \emph{distortion risk measure} (see, e.g.,~\cite{SBR10}) can be seen as a natural extension of the discrete OWA to the continuous case. It can be used to handle more cases in decision making under risk and contains the continuous OWA operator proposed in~\cite{YA04} as a special case.

In the following, we  use the concept of the  \emph{Value at Risk}~(see, e.g.,~\cite{P00}). Let $X$ be a random variable with the cumulative distribution function (CDF) $F_X$. The Value at Risk (VaR for short), for $t \in (0,1]$,  is defined as
\begin{equation}
 \VAR_t(X) = \inf\{y \in \mathbb{R} : F_X(y) \geq t\}= \inf\{y \in \mathbb{R} : \Pr[X\leq y] \geq t\}. 
 \label{evar}
 \end{equation}
We also define $\VAR_0(X) = \inf\{y \in \mathbb{R} : F_X(y) > 0\}$ to handle the boundary case $t=0$. 
The Value at Risk is the quantile function of~$F_X$.  If $X$ has a bounded support, then $\VAR_t(X)$ is bounded for every $t\in [0,1]$.

\subsection{Discrete OWA operator}
\label{sec:discowa}
Let  $\pmb{v}\in\R_{+}^K$ be a vector of reals that can represent, for example, a finite number of realizations of some uncertain quantity or the values of $K$ objective functions that are to be aggregated.
Let $\pmb{w}=(w_1,\dots,w_K)\in \R_{+}^K$ be a vector of nonnegative weights such that  
$w_1+w_2+\dots+w_K=1$. The \emph{Ordered Weighted Averaging} of $\pmb{v}\in \R^K$ is defined as~\cite{YA88}:
\begin{equation}
\OWA^{\mathrm{d}}_{\pmb{w}}(\pmb{v}) =\sum_{i\in [K]} w_i v_{\pi(i)},
\label{owawdef}
\end{equation}
where $[K]=\{1,\ldots,K\}$ and
$\pi$ is a permutation of the set~$[K]$ such that
$v_{\pi(1)}\geq v_{\pi(2)}\geq \dots \geq v_{\pi(K)}$. In this paper, we  call~(\ref{owawdef}) a \emph{discrete OWA operator}.

It is worth pointing out that $\OWA^{\mathrm{d}}_{\pmb{w}}$
encompasses several notable criteria used in decision making as special cases. Namely,
if $w_1=1$ and $w_i=0$ for $i=2,\dots,K$, then $\OWA^{\mathrm{d}}_{\pmb{w}}$ becomes the maximum.
 If $w_K=1$ and $w_i=0$ for $i=1,\dots,K-1$, then $\OWA^{\mathrm{d}}_{\pmb{w}}$ becomes the minimum. In general, if $w_k=1$ and $w_i=0$ for $i\in [K]\setminus\{k\}$, then OWA is the $k$-th largest
component of~$\pmb{v}$.
 In particular, when $k=\lfloor K/2 \rfloor +1$, the $k$th largest component is the median.
If $w_i=\frac{1}{K}$ for all $i \in [K]$, then $\OWA^{\mathrm{d}}_{\pmb{w}}$ is the average. Finally, if $w_1=\lambda$ and $w_K=1-\lambda$, for some fixed $\lambda\in [0,1]$, and $w_i=0$ for the remaining weights, then we get the well-known 
Hurwicz criterion. In practice, the case with nonincreasing weights, $w_1\geq w_2\geq\dots\geq w_K$, is important when minimizing an OWA objective.  Nonincreasing weights can model risk-averse decision makers, who assign larger weights to larger components of $\pmb{v}$. The discrete OWA operator with nonincreasing weights also has some favorable computational properties (see, e.g.,~\cite{OS03, BKKZ25, KZ15}) which we exploit later in this paper.

It was suggested in~\cite{YA04} that the weight vector $\pmb{w}$ in  $\OWA^{\mathrm{d}}_{\pmb{w}}$ can be parameterized by a function
$Q: [0, 1]\rightarrow [0, 1]$ such that $Q(0)= 0$, $Q(1)= 1$, and $Q(x)\geq Q(y)$ for $x>y$, called a \emph{basic unit monotone} (BUM) function. The BUM function $Q$ induces a weight vector $\pmb{w}\in \R^K_{+}$, such that 
\[
w_i=Q\left(\frac{i}{K}\right)-Q\left(\frac{i-1}{K}\right),  i\in [K].
\] 
It is easily seen that $w_i\geq 0$, $i\in [K]$, and $\sum_{i\in[K]} w_i = 1$. The BUM function $Q$ does not need to be continuous. However, it has bounded variation, since it is  bounded and nondecreasing.
The application of the BUM function leads to a version of  the discrete OWA operator~(\ref{owawdef}), parameterized by~$Q$:
\begin{equation}
\label{owadef}
\OWA^{\mathrm{d}}_{Q}(\pmb{v}) = \sum_{i\in[K]} \left(Q\left(\frac{i}{K}\right)-Q\left(\frac{i-1}{K}\right)\right) v_{\pi(i)}.
\end{equation}

The discrete OWA operator can be naturally  formulated  within a stochastic framework. 
Let $\pmb{v}=(v_1,\dots,v_K)\in \R^K$ represent $K$ independent and identically distributed (i.i.d.) samples
of a random variable~$X$ with CDF~$F_{X}$. The empirical distribution function
associated with $\pmb{v}$ is defined as:
\begin{equation}
F_K(x)=\frac{1}{K}|\{ i\,:\, v_i\leq x, i\in [K]\}|.
\label{ecdf}
\end{equation}
An estimator of $\VAR_{1-t}(X)$  is then given by (see, e.g.,~\cite{serfling2009approximation}):
\begin{equation}
 \VAR^{(K)}_{1-t}(X) = F^{-1}_K(1-t)= \inf\{y \in \mathbb{R} : F_K(y) \geq 1-t\}. 
 \label{eevar}
 \end{equation}
Furthermore, 
 if the samples are ordered so that
$v_{\pi(1)}\geq v_{\pi(2)}\geq\cdots\geq v_{\pi(K)}$, then
$\VAR^{(K)}_{1-t}(X)=v_{\pi(\lfloor tK \rfloor+1 )}$ for $t\in [0,1)$ and $\VAR^{(K)}_{1-t}(X)=v_{\pi(K)}$ for $t=1$.
Therefore, the discrete operator $\OWA^{\mathrm{d}}_{Q}$ can be expressed as a weighted sum of the empirical quantiles of~$X$:
\begin{equation}
\label{owaequiv}
   \OWA^{\mathrm{d}}_{Q}(\pmb{v})=\sum_{i\in[K]} \left(Q\left(\frac{i}{K}\right)-Q\left(\frac{i-1}{K}\right)\right) \VAR^{(K)}_{1-\frac{i-1}{K}}(X).
 \end{equation}

In the next subsection, we use
the representation~(\ref{owaequiv}) to extend the OWA operator to the continuous case.

\subsection{Distortion risk measure - an extension of the OWA operator to the continuous case}
\label{ssrm}

The discrete OWA operator, described in Section~\ref{sec:discowa}, aggregates a finite number of components of a given vector~$\pmb{v}$. However, in many applications, we wish to aggregate an infinite number of values contained, for example, in a closed interval, which serves as a simple model of uncertainty~\cite{KY97}.
It turns out that the discrete operator $\OWA^{\mathrm{d}}_{Q}$ can be extended to the continuous case employing the concept of a \emph{distortion risk measure} introduced in~\cite{D90} (see also~\cite{SBR10}). In general, a \emph{risk measure} is a function $\rho$ that assigns a real value to a random variable $X$, which typically denotes a random loss (or cost)~\cite{KAM11}.
From now on, we make the following assumption on $X$:
\begin{enumerate}[label={(}A\arabic*{)}]
\item The random variable $X$ has a bounded support, i.e. there exist $\underline{x},\overline{x}\in \R$ such that
$\Pr[X\in [\underline{x},\overline{x}]]=1$, and
$\VAR_t(X)$ is  continuous on $t\in [0,1]$.
\label{A1}
\end{enumerate}

The following definition, formulated using the notions already introduced, provides a natural extension of~(\ref{owaequiv}) to the continuous case.
\begin{definition}[\cite{D90, SBR10}]
\label{defnewowa}
The distortion risk measure for $X$ is defined as
\begin{equation}
\label{newowa1}
\Sp_{Q}(X) = \int_0^1 \VAR_{1-t}(X)\, d Q(t).
 \end{equation}
If $Q$ is continuously differentiable in $[0,1]$, then 
\begin{equation}
\label{newowa2}
\Sp_{w}(X) = \int_0^1 w(t)\VAR_{1-t}(X)\, d t,
 \end{equation}
where $w(t)=\frac{d Q(t)}{d t}$, $t\in [0,1]$, is the weight function induced by $Q$.
\end{definition}

Note that the Riemann-Stieltjes integral~(\ref{newowa1}) exists
due to assumption~\ref{A1} (see, e.g.,~\cite{R76}). It is easy to check that the weight function $w(t)$ is normalized, that is $\int_{0}^1 w(t) dt=Q(1)-Q(0)=1$. 
The special case $\Sp_{w}(X)$ is called a \emph{spectral risk measure}~\cite{A02,KPB21}. In stochastic terminology, $Q$ is called a \emph{distortion function} and the corresponding weight function $w$ is called a \emph{spectrum}. In this paper, we will follow the terminology from~\cite{YA04} and continue to call $Q$ a BUM function and $w$ a weight function, respectively.

\begin{figure}[ht]
\includegraphics[height=5cm]{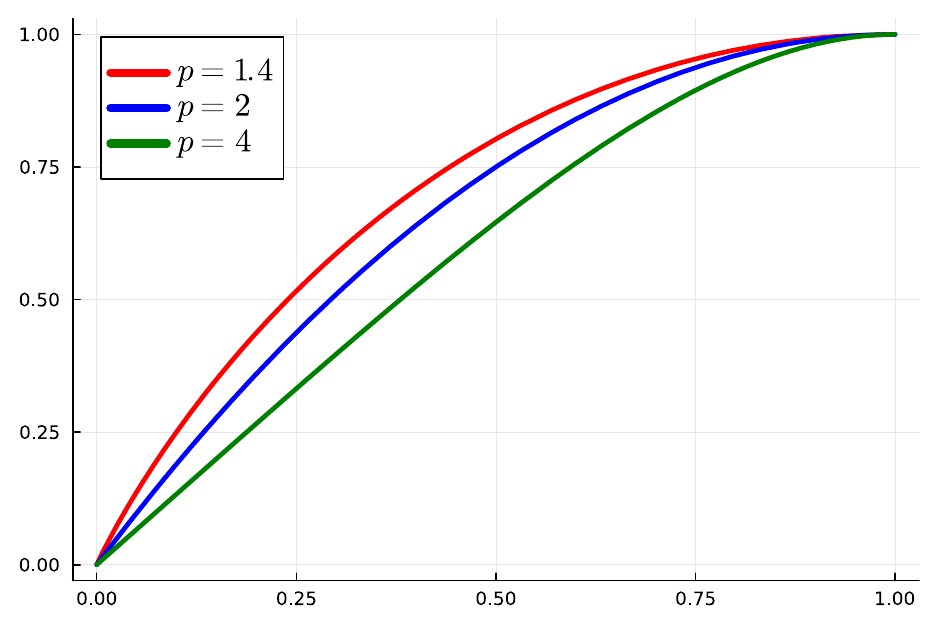}
\includegraphics[height=5cm]{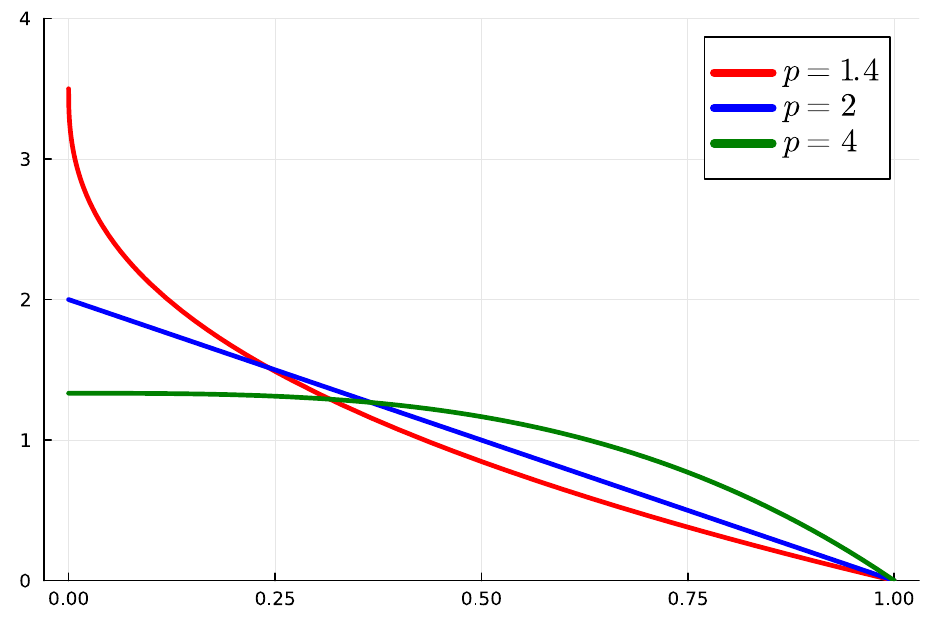}
\caption{The BUM function $Q_p(t)$ and the corresponding nonincreasing weight function $w(t)=\frac{d Q(t)}{dt}$.} \label{figsampq}
\end{figure}

Just like the discrete OWA operator,
the distortion risk measure~$\Sp_{Q}$ generalizes some well-known criteria in stochastic optimization.  If $Q(t)=0$ for $t\in [0,\alpha]$ and $Q(t)=1$ for $t\in(\alpha,1]$ (i.e., $Q$ is a step function), then $\Sp_{Q}(X) = \VAR_{1-\alpha}(X)$. For the boundary cases, we get $\Sp_Q(X)=\underline{x}$ (if $\alpha=1$) and $\Sp_Q(X)=\overline{x}$ (if $\alpha=0$).
From a practical and computational point of view, the case where $Q$ is a concave function is very important. A concave BUM function $Q$ models risk-averse decision makers who overestimate the probabilities of large losses. It represents a level of aversion towards uncertainty (see, e.g.,~\cite{WYP97}).
It is easy to verify that if $Q$ is concave and continuously differentiable, then the corresponding weight function $w$ is nonincreasing. Hence, likewise in the discrete OWA operator with nonincreasing weights, larger weights are assigned to large values of the loss. A commonly used concave BUM function is $Q_\alpha(t)=\min\left\{\frac{t}{\alpha},1\right\}$ for a fixed  $\alpha\in(0,1]$. Then
$$\Sp_{Q}(X) = \frac{1}{\alpha}\int_0^\alpha\VAR_{1-t}(X)\, d t$$
 is the \emph{Conditional Value at Risk}~\cite{rockafellar2000optimization}. In particular, for $\alpha=1$ we get the expectation $\Sp_{Q}(X)=\E[X]$. Consider the following BUM function:
\begin{equation}
Q_p(t)=\frac{p}{p-1}\left(t-\frac{1}{p}t^p\right), \;p> 1.
\label{eqp}
\end{equation} 
It is easy to check that $Q_p(t)$ is concave and continuously differentiable.
The corresponding weight function is $w(t)=\frac{p}{p-1}(1-t^{p-1})$, $t\in [0,1]$ (see Figure~\ref{figsampq}). The smaller $p$, the more risk-averse the decision-maker is. Observe that the weight function $w$ is convex for $p<2$, linear if $p=2$ and concave for $p>2$.  

Suppose that a random variable~$X$ follows a uniform probability distribution on the interval~$[a, b]$, 
which models, for example, a principle of insufficient reason~\cite{LR57}.
Then $\VAR_{1-t}(X)=b-t(b-a)$ for $t\in [0,1]$. Therefore, for a continuously differentiable BUM function $Q$, we get by~(\ref{newowa2}):
\begin{equation}
\label{owayag}
\Sp_{Q}(X)=\int_0^1 \frac{dQ(t)}{dt} (b-t(b-a))\,d t,
\end{equation}
which is exactly  the continuous OWA operator proposed in~\cite{YA04} to aggregate the values in a closed interval~$[a, b]$. However, the application of~(\ref{owayag}) to optimization is limited, because it does not extend to the case where $X$ is not uniformly distributed. 
Such a situation naturally arises in optimization problems with a linear objective function with random coefficients.
In this case, the solution cost is the sum of  random variables, which need not be uniformly distributed even if the cost coefficients themselves are.
We will address this issue later in the paper.

\subsection{Properties of distortion risk measure}
\label{sec:prop}

In this section, we recall known properties of the distortion risk measure and prove some new ones, which will be used later, in the optimization context.

\begin{property}[Boundeness]
\label{propbound}
The distortion risk measure is lower-bounded and upper-bounded,
\begin{equation}
\underline{x} \leq \Sp_Q(X) \leq \overline{x}.
\label{propbound1}
\end{equation}
If the BUM function $Q$ is concave, then 
\begin{equation}
\E[X]\leq \Sp_Q(X).
\label{propbound2}
\end{equation}
Furthermore, if $X$ is nonnegative r.v., then
\begin{equation}
\E[X]\leq \Sp_Q(X)\leq \left(\sup_{t\in (0,1]}\frac{Q(t)}{t}\right)\E[X]
\label{propbound2a}
\end{equation}
and if, additionally, $X$ has a symmetric probability distribution, then 
 \begin{equation}
 \E[X]\leq  \Sp_Q(X)\leq \min\left\{2,  \sup_{t\in (0,1]}\frac{Q(t)}{t}\right\}\E[X]
 \label{propbound3}
 \end{equation}
 (if $Q$ is continuously differentiable, then, in (\ref{propbound2a}) and~(\ref{propbound3}), 
 $\sup_{t\in (0,1]}Q(t)/t=w(0)$, where $w$ is the weight function induced by $Q$).
  \end{property}
  \begin{proof}
  For any BUM function $Q$, we get 
  $\underline{x}=\int_0^1 \underline{x} dQ(t)\leq \int_0^1 \VAR_{1-t}(X) dQ(t)\leq \int_0^1 \overline{x} dQ(t) =\overline{x}$ (see Assumption~\ref{A1}).
  The proof of inequalities~(\ref{propbound2}) and (\ref{propbound2a}) is deferred to Appendix~\ref{dod}.
    Finally, if $X$ has a symmetric probability distribution, then $\E[X]=\frac{1}{2}(\underline{x}+\overline{x})$.
    Since $\Sp_Q(X) \leq \overline{x}$ and $\overline{x}\leq 2\E[X]$,
    $\underline{x}\geq 0$, we have $\Sp_Q(X) \leq 2\E[X]$.
     Combining this with~(\ref{propbound2}) and (\ref{propbound2a}) yields the upper bound~(\ref{propbound3}).
  \end{proof}
\begin{property}
\label{propsym}
If the BUM function $Q$ is such that $Q(t)=1-Q(1-t)$ for every $t\in [0,1]$
and $X$ is a random variable with a symmetric probability distribution, then
\begin{equation}
\Sp_Q(X)=\E[X]. \label{esym}
\end{equation}
\end{property}
\begin{proof}
See Appendix~\ref{dod}.
\end{proof}
Property~\ref{propsym} leads to the following corollary:
\begin{cor}
If  $Q(t)=1-Q(1-t)$ for $t\in [0,1]$ and $Q$ is continuously differentiable, then the weight function~$w$ is symmetric,
$w(t)=w(1-t)$ for $t\in [0,1]$, and $\Sp_Q(X)=\Sp_w(X)=\E[X]$.
\end{cor}
The following property of distortion risk measure is well-known (see, e.g.,~\cite{SBR10}):
\begin{property}
\label{propdist}
The distortion risk measure $\Sp_{Q}$ satisfies the following properties:
\begin{enumerate}
\item (Monotonicity) If $X\leq Y$, then $\Sp_{Q}(X)\leq \Sp_{Q}(Y)$, where $X\leq Y$ means  $\Pr[X\leq z]\geq \Pr[Y\leq z]$ for each $z\in \R$.
\item (Translation invariance) For any constant $a\in \R$, $\Sp_{Q}(X+a)=\Sp_Q(X)+a$.
\item (Positive homogeneity)
For any constant $a\geq 0$, $\Sp_{Q}(aX)=a\cdot\Sp_Q(X)$.
\item (Subadditivity)
If the BUM function $Q$ is concave, then 
	$\Sp_{Q}(X+Y)\leq \Sp_{Q}(X)+\Sp_{Q}(Y)$.

\end{enumerate}
\end{property}
From Property~\ref{propdist} it follows that
$\Sp_Q$  is a \emph{coherent risk measure}, 
provided that the BUM function $Q$ is concave. 
The next theorem shows that $\OWA_{Q}^\textup{d}$ is an estimator of~$\Sp_{Q}$.
\begin{theorem}[Convergence of Discrete OWA]\label{th:convergence} 
Let $\pmb{v}^{(K)}=(v^{(K)}_1,\dots,v^{(K)}_K)$ be  $K$ i.i.d. 
samples
of a random variable~$X$ with the CDF~$F_{X}$. 
Then
\[
\lim_{K \to \infty} \OWA_{Q}^\textup{d}(\pmb{v}^{(K)}) = \Sp_{Q}(X) \quad\text{(almost surely)}.
\]
\end{theorem}
\begin{proof}
We recall that the empirical distribution function $F_K(x)$ associated with $\pmb{v}^{(K)}$ is defined by~(\ref{ecdf}), 
and that $\VAR^{(K)}_{1-t}(X)$ and $F^{-1}_X(x)$ are defined by~(\ref{eevar}).
By~(\ref{owadef}), (\ref{owaequiv}), the definition of the weights
$w_i = Q(i/K) - Q((i-1)/K)$,   $i \in [K]$,
and using the fact that $v^{(K)}_{\pi(i)} = \VAR^{(K)}_{1 - \frac{i-1}{K}}(X)$, we obtain
\begin{align*}
&  |\OWA_{Q}^\textup{d}(\pmb{v}^{(K)}) - \Sp_{Q}(X)|= | \sum_{i\in[K]} w_i v^{(K)}_{\pi(i)} - \int_0^1\VAR_{1-t}(X) \, dQ(t)| \\
&= | \sum_{i\in[K]} w_i v^{(K)}_{\pi(i)} - \int_0^1 \VAR_{1-t}(X) \, d Q(t) \\
&\qquad + \sum_{i\in[K]} w_i \VAR_{1-\frac{i-1}{K}}(X) - \sum_{i\in[K]} w_i \VAR_{1-\frac{i-1}{K}}(X)  | \\
&\leq \underbrace{| \sum_{i\in[K]} w_i \VAR^{(K)}_{1-\frac{i-1}{K}}(X) - \sum_{i\in[K]} w_i \VAR_{1-\frac{i-1}{K}}(X) |}_{=:f^{(K)}} \\
&\qquad + \underbrace{|\sum_{i\in[K]} (Q(i/K)-Q((i-1)/K) \VAR_{1-\frac{i-1}{K}}(X) - \int_0^1 \VAR_{1-t}(X) \, dQ(t) |}_{=:g^{(K)}}.
\end{align*}
By assumption~\ref{A1},
$\VAR_{1-t}(X)$ is continuous on $t\in[0,1]$. In addition, $Q(t)$ has a bounded variation. Therefore,
the function $\VAR_{1-t}(X)$ is Riemann–Stieltjes integrable with respect to $Q(t)$.
By the definition of the Riemann–Stieltjes integral, we have
$\lim_{K\to\infty} g^{(K)} = 0$.

We now show that $\lim_{K\to\infty} f^{(K)} = 0$ almost surely.
By the triangle inequality and the fact that $\sum_{i\in [K]} w_i = 1$, we obtain
\[
f^{(K)} \leq \sup_{t\in[0,1]}
|\VAR_{1-t}^{(K)}(X) - \VAR_{1-t}(X)|.
\]
Since $\VAR_{1-t}(X)$ is continuous on the interval [0,1], it is uniformly continuous.
By the Glivenko–Cantelli theorem, the empirical distribution function satisfies
$\lim_{K\to\infty} F_K = F_X$ uniformly almost surely (see, e.g.,~\cite{serfling2009approximation}).
This implies a uniform almost sure convergence of empirical quantiles:
\[
\lim_{K\to\infty}\sup_{t\in[0,1]}
| \VAR_{1-t}^{(K)}(X)-\VAR_{1-t}(X)|= 0
\quad\text{(almost surely)}.
\]
Consequently, $\lim_{K\to\infty} f^{(K)} = 0$ almost surely.
Combining both parts, we obtain
\[
\lim_{K\to\infty}
| \OWA_{Q}^\mathrm{d}(\pmb{v}^{(K)}) - \Sp_Q(X)|
= 0
\quad\text{(almost surely)},
\]
which completes the proof.
\end{proof}

%
%

\section{Optimization under risk}
\label{sec:2}

Consider the following generic optimization problem with a linear objective function:
\begin{equation}
\label{defopt}
\begin{array}{lll}
		\min &  \pmb{c}^\top \pmb{x} \\
			&\pmb{x}\in \X\subseteq \R^n
			\end{array}
			\end{equation}
where $\X$ is a set of feasible solutions and $\pmb{c}\in \mathbb{R}^n$ is a  cost vector. In many applications, the cost vector $\pmb{c}$ is uncertain. Under the interval uncertainty representation, we only know that $c_i\in[\underline{c}_i,\overline{c}_i]$ for each $i\in [n]$
(see, e.g.,~\cite{BN09,KY97}).
 Then  $\cU = \prod_{i\in[n]} [\underline{c}_i,\overline{c}_i]$ is an \emph{interval uncertainty set}. Notice that the cost of $\pmb{x}\in \X$ belongs to the interval $[\underline{\pmb{c}}^\top \pmb{x},\overline{\pmb{c}}^\top \pmb{x}]=[\sum_{i\in [n]} \underline{c}_ix_i, \sum_{i\in [n]} \overline{c}_ix_i]$. Assuming that this cost is uniformly distributed, we can apply the continuous OWA operator~(\ref{owayag}), proposed in~\cite{YA04}, to aggregate the values in $[\underline{\pmb{c}}^\top \pmb{x},\overline{\pmb{c}}^\top \pmb{x}]$, which leads to the following criterion:
\begin{equation}
\label{owahur0}
\int_0^1 \frac{dQ(t)}{dt} (\overline{\pmb{c}}^\top \pmb{x}-t(\overline{\pmb{c}}^\top \pmb{x}-\underline{\pmb{c}}^\top \pmb{x}))\;dt = \lambda\overline{\pmb{c}}^\top \pmb{x}+(1-\lambda)\underline{\pmb{c}}^\top \pmb{x},
\end{equation}
where $\lambda = 1- \int_0^1 t w(t) dt$ and $w(t)=\frac{d Q(t)}{dt}$.
Notice that $\lambda\in [0,1]$
and~(\ref{owahur0}) is the well-known Hurwicz criterion (see, e.g.,~\cite{LR57}), which is a convex combination of the minimum and the maximum cost of $\pmb{x}$ in $\mathcal{U}$. Therefore, this type of OWA aggregation is always given by the weighted average of the minimum and maximum possible cost of $\pmb{x}$, defined by the bounds of the cost intervals.  From the computational complexity point of view, it is an advantage that we thus obtain a nominal optimization problem of the form~(\ref{defopt}), where the costs are  $c_i=\lambda \overline{c}_i+(1-\lambda)\underline{c}_i$, $i\in [n]$.
 However, we note two disadvantages: first, the actual shape of the function $Q$ is irrelevant, as it is replaced by a single value $\lambda$. This limits the flexibility to express the decision maker's preferences. Secondly, and more subtly, it is implicitly assumed that all cost outcomes in the interval $[\underline{\pmb{c}}^\top \pmb{x},\overline{\pmb{c}}^\top \pmb{x}]$ are equally likely. This assumption is not true if the cost of $\pmb{x}$ is a sum of random variables, in which case the probability distribution of the cost in $[\underline{\pmb{c}}^\top \pmb{x},\overline{\pmb{c}}^\top \pmb{x}]$ is typically not uniform.

\subsection{Optimization problem with a distortion risk measure}

Let  $\pmb{C}=(C_1,\dots,C_n)^\top$ be a nonnegative random cost vector,
where all its 
components $C_i$, $i\in [n]$, have bounded supports $[\underline{c}_i,\overline{c}_i]$, $i\in [n]$. In what follows, $\pmb{C}$ has a bounded support $\cU = \prod_{i\in[n]} [\underline{c}_i,\overline{c}_i]\subset \R_{+}$, which is an interval uncertainty set. 
In this paper, we address the following optimization problem:
\begin{equation}
\label{defowaprob}
\begin{array}{lll}
		\min &   \Sp_{Q}(\pmb{C}^\top \pmb{x}) \\
			&\pmb{x}\in \X\subseteq \R^n_{+}
			\end{array}
\end{equation}
From now on, we assume that the random variable~$\pmb{C}^\top \pmb{x}$
satisfies~\ref{A1}.
If the continuous normalized weight function $w:[0,1]\rightarrow [0,1]$ ($\int_{0}^1 w(t) dt=1$) is explicitly given or is obtained from a continuously differentiable BUM function~$Q$, then we replace $ \Sp_{Q}(\pmb{C}^\top \pmb{x})$ with $\Sp_{w}(\pmb{C}^\top \pmb{x})$ in~(\ref{defowaprob}). 

\subsection{Solution algorithms}
\label{sec:solalg}

In this section, we will discuss some methods for solving the problem~(\ref{defowaprob}). It is generally computationally intractable, as it is $\#$P-hard to compute the value of $\Sp_{Q}(\pmb{C}^\top \pmb{x})$ for a given solution $\pmb{x}$ (see Corollary~\ref{corcompl2}, later in this paper). Recall that a problem is $\#$P-hard if every problem in the class $\#$P is polynomial-time Turing reducible to it, where the class $\#$P consists of all counting problems associated with decision problems in NP (see~\cite{V79}). The $\#$P-hard problems are at least as hard as NP-complete problems. Therefore, some approximation methods for solving~(\ref{defowaprob}) are required. 
Let us define two special cost vectors $\pmb{c}^{Q}$ and $\hat{\pmb{c}}$, whose components are
\begin{align}
c_i^Q&=\Sp_{Q}(C_i), &i\in [n], \label{cqdef}\\
\hat{c}_i&=\E[C_i], &i\in [n].
\end{align}
Let us also denote $\eta=\sup_{t\in (0,1]}\frac{Q(t)}{t}$. If $Q$ is concave and continuously differentiable, then $\eta=w(0)$.
\begin{lemma}
\label{lem2appr}
If the BUM function $Q$ is concave and  $\pmb{x}^*\in \X\subseteq \R_{+}^n$ is an optimal solution to~(\ref{defopt}) for the cost vector $\pmb{c}^Q$, then
\begin{equation}
\forall \pmb{x}\in \X\;\; \Sp_Q(\pmb{C}^\top\pmb{x}^*)\leq \eta\cdot \Sp_Q(\pmb{C}^\top\pmb{x}),
\label{apsup}
\end{equation}
if additionally $C_i$, $i\in [n]$, have symmetric probability distributions, then
\begin{equation}
\forall \pmb{x}\in \X\;\; \Sp_Q(\pmb{C}^\top\pmb{x}^*)\leq \min\left\{2, \eta\right\}\Sp_Q(\pmb{C}^\top\pmb{x}).
\label{ap2sup}
\end{equation}
\end{lemma}
\begin{proof}
We show only the inequality~(\ref{ap2sup}). The proof of~(\ref{apsup}) is similar. Let $\beta= \min\left\{2,  \eta\right\}$.
By  Property~\ref{propdist} (positive homogeneity and  subadditivity) and Property~\ref{propbound}, we get
$$\Sp_Q(\pmb{C}^\top\pmb{x}^*)\leq \sum_{i\in [n]} x_i^*\Sp_Q(C_i)\leq \sum_{i\in [n]} x_i\Sp_Q(C_i)
\leq \beta\sum_{i\in [n]} x_i \E[C_i]\leq \beta\E[\pmb{C}^\top\pmb{x}]\leq \beta\cdot \Sp_{Q}(\pmb{C}^\top \pmb{x}).$$
\end{proof}
From Lemma~\ref{lem2appr}, we immediately get the following approximation results:
\begin{theorem}
\label{thm2appr}
If the BUM function $Q$ is concave and an optimal solution to the problem~(\ref{defopt}) for $\pmb{c}^Q$
can be computed in  polynomial time,
then the problem~(\ref{defowaprob}) is approximable within $\eta$.
If additionally $C_i$, $i\in [n]$, have symmetric probability distributions, then~(\ref{defowaprob}) is 
  approximable within $\min\{2,\eta\}$.
\end{theorem}
We obtain the same bound as in Theorem~\ref{thm2appr} by using the cost vector~$\hat{\pmb{c}}$  instead of~$\pmb{c}^Q$. Indeed, if $\pmb{x}^*$ is an optimal solution to~(\ref{defopt}) for the cost vector $\hat{\pmb{c}}$, then by the same reasoning as in Lemma~\ref{lem2appr} we get $\Sp_Q(\pmb{C}^\top\pmb{x}^*)
\leq \beta\E[\pmb{C}^\top\pmb{x}^*]\leq \beta\E[\pmb{C}^\top\pmb{x}]\leq \beta\cdot \Sp_{Q}(\pmb{C}^\top \pmb{x})$ for any $\pmb{x}\in \X$ and we note that $\E[\pmb{C}^\top\pmb{x}]=\hat{\pmb{c}}^\top \pmb{x}$. Using $\hat{\pmb{c}}$ can be advantageous, because the expected values of $C_i$, $i\in [n]$, are usually readily available. On the other hand, using $\pmb{c}^Q$ takes into account  the shape of~$Q$ which may yield a tighter bound for particular instances. However, the computation of $\pmb{c}^Q$ may require complex integration.



The next lemma allows us to reduce solving the problem~(\ref{defowaprob}) to solving the underlying deterministic problem~(\ref{defopt}) with the cost vector $\hat{\pmb{c}}$, under an additional assumption on $Q$.
\begin{lemma}
\label{lem:symmetric}
Assume that $C_i$, $i\in [n]$,
are  independent random variables with symmetric probability distributions. 
If the BUM function $Q$ is such that $Q(t)=1-Q(1-t)$ for $t\in [0,1]$, then
\[
\forall\pmb{x}\in\X \quad \Sp_Q(\pmb{C}^\top\pmb{x}) =\E[\pmb{C}^\top\pmb{x}]=\sum_{i\in[n]}\E[C_i]x_i=\hat{\pmb{c}}^\top \pmb{x}.
\]
\end{lemma}
\begin{proof}
The random variable $\pmb{C}^\top\pmb{x}$  is  the sum of $n$~independent, symmetrically distributed r.v's $x_iC_i$, $i\in [n]$. So, it is
symmetrically distributed  around~$\E[\pmb{C}^\top\pmb{x}]=\sum_{i\in[n]}\E[C_i]x_i$ (see, e.g.,~\cite{GS01}).
Thus,
Property~\ref{propsym} yields $\Sp_Q(\pmb{C}^\top\pmb{x}) =\E[\pmb{C}^\top\pmb{x}]=\hat{\pmb{c}}^\top \pmb{x}$.
\end{proof}

Lemma~\ref{lem:symmetric} implies the following result.
\begin{cor}
Assume that the cost vector $\hat{\pmb{c}}$ can be computed in polynomial time and the assumptions of Lemma~\ref{lem:symmetric} are satisfied. Then~(\ref{defowaprob})  can be solved by solving~(\ref{defopt}) for~$\hat{\pmb{c}}$. Therefore, in this case~(\ref{defowaprob}) has the same computational complexity as~(\ref{defopt}).
\end{cor}
\begin{theorem}
If the BUM function $Q$ is concave and $\X$ is a convex set, then $\Sp_{Q}(\pmb{C}^\top\pmb{x})$ is a convex function in $\pmb{x}\in\X$.
\end{theorem}
\begin{proof}
Choose $\lambda\in [0,1]$ and any solutions $\pmb{x}^{(1)},\pmb{x}^{(2)}\in \X$. Let $\pmb{x}=\lambda\pmb{x}^{(1)}+(1-\lambda)\pmb{x}^{(2)}$. By the convexity of~$\X$, $\pmb{x}\in \X$ and by the linearity of $\pmb{C}^\top\pmb{x}$, we have
$\pmb{C}^\top\pmb{x}=\lambda \pmb{C}^\top\pmb{x}^{(1)}+(1-\lambda)\pmb{C}^\top\pmb{x}^{(2)}$, where
$\pmb{C}^\top\pmb{x}^{(1)}$ and $\pmb{C}^\top\pmb{x}^{(2)}$ are two r.v's.
Hence and from Property~\ref{propdist}, we obtain
$$
\Sp_{Q}(\pmb{C}^\top\pmb{x})=\Sp_{Q}(\lambda \pmb{C}^\top\pmb{x}^{(1)}+(1-\lambda)\pmb{C}^\top\pmb{x}^{(2)})
\leq \lambda \cdot \Sp_{Q}(\pmb{C}^\top\pmb{x}^{(1)})+(1-\lambda)\cdot \Sp_{Q}(\pmb{C}^\top\pmb{x}^{(2)}),
$$
which completes the proof.
\end{proof}
\begin{cor}
\label{corconv}
If the BUM function $Q$ is concave and the set  $\X$ is convex, then~(\ref{defowaprob}) is a convex optimization problem.
\end{cor}

Let us now focus on computing the value of
$\Sp_{Q}(\pmb{C}^\top \pmb{x})$ for a given $\pmb{x}\in\X$.
This value can be approximated using the discrete OWA estimator~(\ref{owaequiv}).
Let $(\pmb{c}_i)_{i\in [K]}$ be $K$ i.i.d. samples of the random cost vector
$\pmb{C}$ with the distribution function~$F_{\pmb{C}}$.
They  induce $K$ i.i.d. samples
$(\pmb{c}_i^\top \pmb{x})_{i\in[K]}$ of the random variable
$\pmb{C}^\top \pmb{x}$.
If the samples are ordered so that
$\pmb{c}_{\pi(1)}^\top \pmb{x} \geq \pmb{c}_{\pi(2)}^\top \pmb{x} \geq \cdots
\geq \pmb{c}_{\pi(K)}^\top \pmb{x}$,
then (\ref{owaequiv}) can be written as:
\begin{align}
\OWA_{Q}^\textup{d}((\pmb{c}_i^\top \pmb{x})_{i\in[K]})
&= \sum_{i=1}^K ( Q(i/K) - Q((i-1)/K) )
\VAR^{(K)}_{1-\frac{i-1}{K}}(\pmb{C}^\top\pmb{x}) \label{lcowa} \\
&= \sum_{i=1}^K ( Q(i/K) - Q((i-1)/K))
\pmb{c}_{\pi(i)}^\top \pmb{x}, \label{lcowaw}
\end{align}
where $\VAR^{(K)}_{1-t}(\pmb{C}^\top\pmb{x})$ is an estimator of~$\VAR_{1-t}(\pmb{C}^\top\pmb{x})$,
based on the empirical distribution associated with~$\pmb{C}^\top\pmb{x}$. Observe that $\VAR^{(K)}_{1-t}(\pmb{C}^\top\pmb{x})=\pmb{c}_{\pi(i)}^\top\pmb{x}$ for each $t\in [(i-1)/K, i/K)$ and $i\in [K]$.
Hence, the formula~(\ref{lcowa}) is  the left 
endpoint rule  for approximating the  Riemann-Stieltjes integral  representation of~$\Sp_Q$. Because $\VAR^{(K)}_{1-\frac{2i-1}{2K}}(\pmb{C}^\top\pmb{x})=\pmb{c}_{\pi(i)}^\top\pmb{x}$, $i\in[K]$, we get exactly the same formula as~(\ref{lcowaw}) by applying the midpoint rule:
\begin{align}
   \OWA^{\mathrm{d}}_{Q}((\pmb{c}_i^\top\pmb{x})_{i\in[K]})&=
   \sum_{i\in[K]} (Q(i/K)-Q((i-1)/K)) \VAR^{(K)}_{1-\frac{2i-1}{2K}}(\pmb{C}^\top\pmb{x}) \label{mcowa}\\
   &=\sum_{i\in[K]} (Q(i/K)-Q((i-1)/K)) \pmb{c}_{\pi(i)}^\top\pmb{x}. \label{mcowaw}
\end{align}
By Theorem~\ref{th:convergence}, both 
(\ref{lcowa}) and (\ref{mcowa}) converge almost surely to
$\Sp_Q(\pmb{C}^\top\pmb{x})$.
Moreover, the midpoint rule~(\ref{mcowa}) achieves an approximation error that is not worse, and typically smaller, 
than the error of the trapezoidal rule.
If the BUM function $Q$ is continuously differentiable in $[0,1]$ and
$w(t)=\frac{dQ(t)}{dt}$ for $t\in[0,1]$, then
the midpoint rule approximation of $\Sp_w(\pmb{C}^\top\pmb{x})$ takes the form of
\begin{equation}
   \OWA^{\mathrm{d}}_{w}((\pmb{c}_i^\top\pmb{x})_{i\in[K]})=\sum_{i\in[K]} \frac{1}{K} w\left(\frac{2i-1}{2K}\right)  \pmb{c}_{\pi(i)}^\top\pmb{x}.
    \label{mcowaww}
\end{equation}
In \cite{KPB21}, the case where a weight function~$w$ is available was analyzed, and the trapezoidal rule was used to approximate~$\Sp_w$, together with concentration bounds derived under smoothness assumptions on both~$w$ and $\VAR_{t}$.
It is well known that the midpoint rule has an approximation error of the same asymptotic order $O(1/K^{2})$ as the trapezoidal rule, but with  half of the error constant, making the midpoint rule  strictly more accurate in this setting.

We are  now ready to provide an algorithm for solving~(\ref{defowaprob}) with the approximated objective function
$\Sp_Q(\pmb{C}^\top\pmb{x})$. We generate $K$ i.i.d. samples $(\pmb{c}_i)_{i\in [K]}$ of the random cost vector
$\pmb{C}$  
and replace $\Sp_Q(\pmb{C}^\top\pmb{x})$ in~(\ref{defowaprob})
with $\OWA^{\mathrm{d}}_{Q}((\pmb{c}_i^\top\pmb{x})_{i\in[K]})$ (see~(\ref{mcowaw}))
or with $\OWA^{\mathrm{d}}_{w}((\pmb{c}_i^\top\pmb{x})_{i\in[K]})$ (see~(\ref{mcowaww})).
This sampling algorithm is shown in the form of Algorithm~\ref{alg:sampling}.
\begin{algorithm}
\caption{Sampling Algorithm (SA)}\label{alg:sampling}
Generate  $K$ i.i.d. samples $(\pmb{c}_i)_{i\in [K]}$ of~$\pmb{C}$\;
$\displaystyle \pmb{x}^* \gets   \operatorname*{arg}\,\min_{\pmb{x} \in \X}
\OWA^{\mathrm{d}}_{Q}((\pmb{c}_i^\top\pmb{x})_{i\in[K]})$
\tcp{or $\displaystyle \operatorname*{arg}\,\min_{\pmb{x} \in \X}
\OWA^{\mathrm{d}}_{w}((\pmb{c}_i^\top\pmb{x})_{i\in[K]})$}\label{l2}
\Return{$\pmb{x}^*$}\;
\end{algorithm}

It remains to explain how to solve the corresponding discrete OWA problem in line~\ref{l2} of Algorithm~\ref{alg:sampling}. 
For this purpose, several methods have already been described in the literature. In the general case, we can use a mixed integer programming reformulation of the discrete OWA \cite{OO12, OS03}. However, if the BUM function $Q$ is concave (the corresponding weights $\pmb{w}$ are nonincreasing), we can use a linear programming reformulation proposed in~\cite{CG15}. Namely, the problem in line~~\ref{l2}  can be formulated as the following optimization problem:
\begin{align*}
\min\ & \sum_{k\in[K]} \alpha_k + \beta_k \\
\text{s.t. } & \alpha_k + \beta_j \ge w_k\pmb{c}_j^\top\pmb{x} & \forall j,k\in[K] \\
& \pmb{x}\in\X \\
& \pmb{\alpha},\pmb{\beta}\in\mathbb{R}^K
\end{align*}
where $w_k=Q(k/K)-Q((k-1)/K)$
(or $w_k=\frac{1}{K} w\left(\frac{2k-1}{2K}\right)$)  for $k\in [K]$.
If $\X$ is a convex set, then we get a convex optimization problem. In particular, if $\X$ is described by linear constraints, we get a linear programming problem. In this case, we can use a large sample size $K$ in Algorithm~\ref{alg:sampling}.

\subsection{The problem with independent uniformly distributed cost coefficients}

\subsubsection{Motivation}

In this section, we investigate a special, relevant case of problem~(\ref{defowaprob}), where the components of  $\pmb{C}=(C_1,\dots,C_n)^\top$ are independent and the random cost $C_i$ is uniformly 
distributed in $[\underline{c}_i,\overline{c}_i]$ where $\underline{c}_i\geq 0$ for each $i\in [n]$.
In this case, the cost vector $\pmb{C}$ is uniformly distributed over the uncertainty set 
$\mathcal{U} = \prod_{i\in [n]} [\underline{c}_i, \overline{c}_i]\subset \R^n_{+}$. 
This corresponds to the situation in which we assume complete ignorance about the values of the cost coefficients. Then the principle of insufficient reason leads to uniform distribution for the uncertain costs (see, e.g.,~\cite{LR57}).
However, note that the random variable $\pmb{C}^\top \pmb{x}$  is
 not uniformly distributed.
The CDF of a sum of $n$ uniformly distributed random variables is a (scaled) version of the Irwin–Hall distribution, which is a continuous, increasing, piecewise polynomial of degree~$n$ (see~\cite[Chapter~9]{johnson1995continuous}). Hence, applying the criterion~(\ref{owahur0}), proposed in~\cite{YA04}, can be inappropriate in this case and using  $\Sp_Q(\pmb{C}^\top \pmb{x})$ can better reflect the uncertainty. Notice that the random variable 
 $\pmb{C}^\top \pmb{x}$ has 
 the support $[\underline{\pmb{c}}^\top \pmb{x},\overline{\pmb{c}}^\top \pmb{x}]\subset \R_+$ and thus satisfies the assumption~\ref{A1}.
 
\begin{figure}[H]
    \centering
    \includegraphics[width=0.8\textwidth]{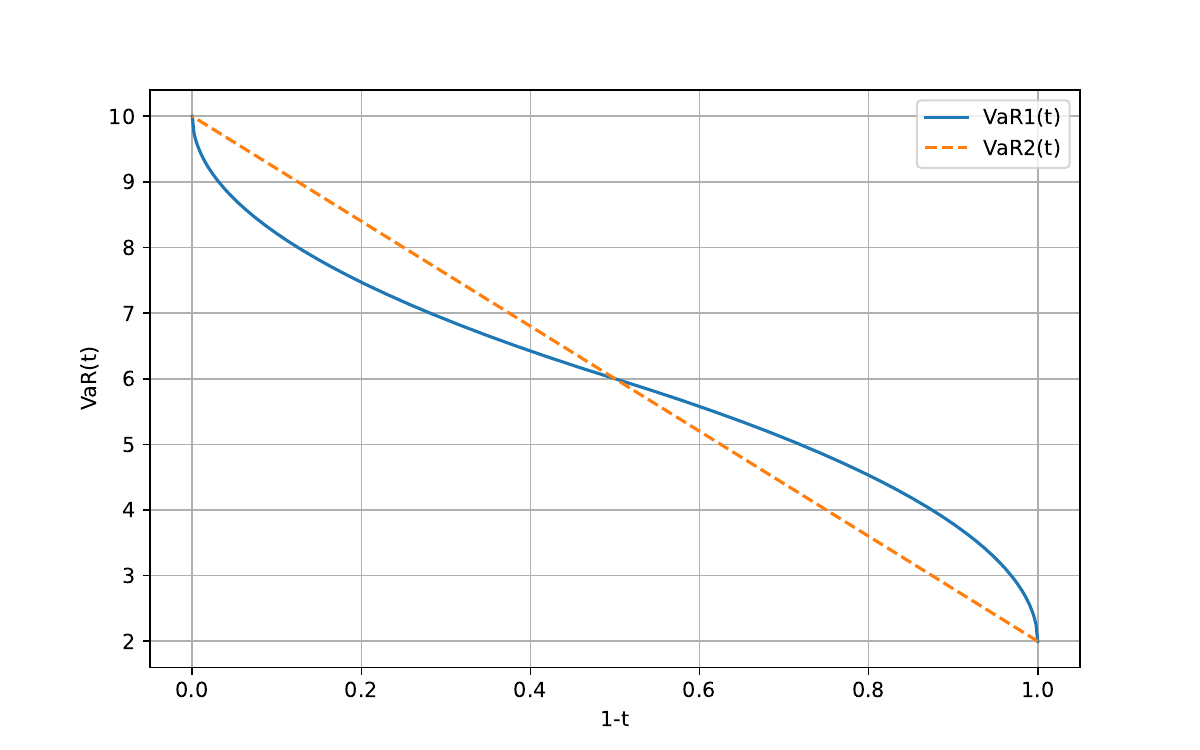}
    \caption{Comparison of VaR profiles for different solutions.}
    \label{fig:VaR}
    \end{figure}

Let us illustrate the model by a simple example.
Consider the optimization problem~(\ref{defopt}), where $\X = \{\pmb{x}\in\{0,1\}^3 : x_1 - x_2 = 0,\ x_1 + x_3 = 1\}$. The set $\mathcal{X}$ represents a situation in which we have three items and can choose both items~1 and 2, or only item~3. Let 
$C_1\sim U[1,5]$, $C_2\sim U[1,5]$, and $C_3\sim U[2,10]$ so $\pmb{C}^\top
=(C_1,C_2,C_3)$ is uniformly distributed in $\mathcal{U}=[1,5]\times [1,5]\times [2,10]$.
Choose the solution $\pmb{x}^{(1)}=(1,1,0)$ in which we select items~1 and~2. The CDF of the resulting cost distribution is then:
\[ \Pr[\pmb{C}^\top\pmb{x}^{(1)} \le y] = \begin{cases}
    0 & \text{ if } y \le 2 \\
    \frac{(y-2)^2}{32} & \text{ if } y\in [2,6] \\
    1-\frac{(10-y)^2} {32} & \text{ if } y \in[6,10] \\
    1 & \text{ else }
    \end{cases}\]
The Value at Risk is the inverse of the CDF. That is,
\[ \VAR_t(\pmb{C}^\top\pmb{x}^{(1)}) = \begin{cases}
    \sqrt{32t} + 2 & \text{if } t \in [0, 0.5] \\
    10 - 4 \sqrt{2 - 2t} & \text{if } t \in [0.5, 1] \\
    \end{cases} \] 
Choose now  solution $\pmb{x}^{(2)}=(0,0,1)$ in which we choose item~3. The CDF of the resulting cost distribution is then
\[  \Pr[\pmb{C}^\top\pmb{x}^{(2)} \le y] = \begin{cases} 0 & \text{ if } y < 2 \\
\frac{(y-2)}{8} & \text{ if } y\in [2,10] \\
1 & \text{ if } y > 10 \end{cases}\]
and the Value at Risk is given by
\[ \VAR_t(\pmb{C}^\top\pmb{x}^{(2)}) = 2 + 8t \]
for $t \in [0, 1]$.
Figure~\ref{fig:VaR} shows the Value at Risk curves for the two solutions in our problem. VaR1($t$) represents the strategy of selecting items 1 and 2, while VaR2($t$) represents the selection of only item 3. The curves intersect at the extremes, corresponding to the best case and worst case scenarios, respectively, and also in the middle, which represents the average case. However, solution $\pmb{x}^{(1)}$ is generally more favorable for risk-averse decision-makers who put more weight on worse outcomes.
This example illustrates a difference between solutions $\pmb{x}^{(1)}$ and $\pmb{x}^{(2)}$ that is lost if we consider these solutions only through the lens of the lower bound $a=2$ and upper bound $b=10$ that they both achieve. If we use a risk-neutral weight function $w(t)=1$, then
\[ \Sp_{w}(\pmb{C}^\top\pmb{x}^{(1)}) = \Sp_{w}(\pmb{C}^\top\pmb{x}^{(2)}) = 6. \]
However, if we emphasize results with high (bad) objective values more than results with low (good) objective values, the picture becomes different. For example, using the more risk-averse weight function $w(t)=3(1-t)^2$, we obtain
\[ \Sp_{w}(\pmb{C}^\top\pmb{x}^{(1)}) = 7.4 \text{ and } \Sp_{w}(\pmb{C}^\top\pmb{x}^{(2)}) = 8 \]
and we recognize that the solution $\pmb{x}^{(1)}$ performs better than the solution $\pmb{x}^{(2)}$.


\subsubsection{The computational complexity of the problem}
\label{sec:3}

We now characterize the computational complexity of the problem~(\ref{defowaprob}) with the independent, uniformly 
distributed objective coefficients. We first show that the problem of computing the value of 
$\Sp_{Q}(\pmb{C}^\top\pmb{x})$ for a given solution~$\pmb{x}$ is $\#$P-hard.
We will use the following result~\cite[Remark, p. 970]{DF88}.
\begin{lemma}
\label{corcompl}
Let $\pmb{a}\in \N^n$ and $b\in \N$.
Suppose that $\pmb{Z}=(Z_1,\ldots,Z_n)$ is
a random vector whose components~$Z_i$, $i\in[n]$, are independent and uniformly distributed on~$[0,1]$
Then computing the value $\Pr[\pmb{a}^\top\pmb{Z}\leq b]$ is $\#$P-hard.
 \end{lemma}
\begin{lemma}
\label{thmcompl}
	Computing the value of $\VAR_t(\pmb{C}^\top \pmb{x})$, for a given solution $\pmb{x}\in \R^n_{+}$ and arbitrary $t\in [0,1]$, is $\#$P-hard.
\end{lemma}
\begin{proof}
Consider the instance from Lemma~\ref{corcompl}: 
we are given $\pmb{a}\in\N^n$ and $b\in\N$, and 
$\pmb{Z}=(Z_1,\dots,Z_n)$ whose components are independent and uniformly distributed over $[0,1]$.
Define $C_i := a_i Z_i$ for $i\in[n]$, so that the random vector 
$\pmb{C}=(C_1,\dots,C_n)^\top$ has independent components uniformly distributed on $[0,a_i]$.
Fix $\hat{\pmb{x}}=(1,1,\dots,1)$.  
Then for every $y\ge 0$ we have
\[
\Pr[\pmb{C}^\top \hat{\pmb{x}} \leq y]
= \Pr\!\left[\sum_{i=1}^n a_i Z_i \leq y\right]
= \Pr[\pmb{a}^\top \pmb{Z} \leq y].
\]
We first show that the following decision problem is $\#$P-hard:
\begin{center}
\emph{Given $t\in[0,1]$ and $b\in\R$, decide whether $\VAR_t(\pmb{C}^\top \hat{\pmb{x}})\leq b$.}
\end{center}
Assume, for the sake of contradiction, that we have an oracle 
\textsc{D-VaR}$(\pmb{C}^\top\hat{\pmb{x}},t,b)$
that answers the above question in polynomial time.
Since $\pmb{a}\in\N^n$ and $b\in\N$, it is known (see, e.g.,~\cite{DF88})
that
\[
\Pr[\pmb{a}^\top \pmb{Z} \leq b] = \frac{p}{q}
\qquad\text{for some integers }p,q,
\]
where we can take  $q = n! \prod_{i=1}^n a_i$.
In particular, every candidate value of $\Pr[\pmb{a}^\top \pmb{Z} \leq b]$ is of the form 
$p/q$ with $p\in\{0,1,\dots,q\}$.
Moreover,
\[
\Pr[\pmb{a}^\top \pmb{Z} \leq b] \geq t
\quad\Longleftrightarrow\quad 
\VAR_t(\pmb{a}^\top\pmb{Z}) \leq b.
\]
We can therefore compute the numerator $p$ using binary search over the set 
$\{0,1,\dots,q\}$, performing $O(\log q)$ oracle calls.
In each step of the binary search, we query the oracle with 
$t = \ell/q$, where $\ell$ is the current midpoint.
Hence, if \textsc{D-VaR} runs in polynomial time, then 
$\Pr[\pmb{a}^\top \pmb{Z} \le b]$ can also be computed in polynomial time.
This contradicts the $\#$P-hardness of computing the value 
$\Pr[\pmb{a}^\top \pmb{Z} \leq b]$ established in Lemma~\ref{corcompl}.
Therefore, the decision version of VaR is $\#$P-hard and in consequence computing
the value of $\VAR_t(\pmb{C}^\top \pmb{x})$ is $\#$P-hard as well.
%
%
%
\end{proof}
\begin{cor}
	\label{corcompl1}
	If the costs $C_i$, $i\in [n]$, are  independent and have uniform probability distributions, then computing the value of $\Sp_{Q}(\pmb{C}^\top \pmb{x})$ for a given solution $\pmb{x}$ is $\#$P-hard.
\end{cor}
\begin{proof}
It is easy to check that
the Value at Risk is a special case of $\Sp_Q(\pmb{C}^\top \pmb{x})$. Indeed,
$\Sp_Q(\pmb{C}^\top \pmb{x})=\VAR_{1-\alpha}(\pmb{C}^\top \pmb{x})$ for the BUM function defined by
$Q(t)=0$ for $t\in [0,\alpha)$ and $Q(t)=1$ for $t\in [\alpha,1]$.
\end{proof}
\begin{cor}
	\label{corcompl2}
	The problem~(\ref{defowaprob}) with the independent uniformly 
distributed objective coefficients  is $\#$P-hard.
\end{cor}

\begin{figure}[ht]
\includegraphics[height=5cm]{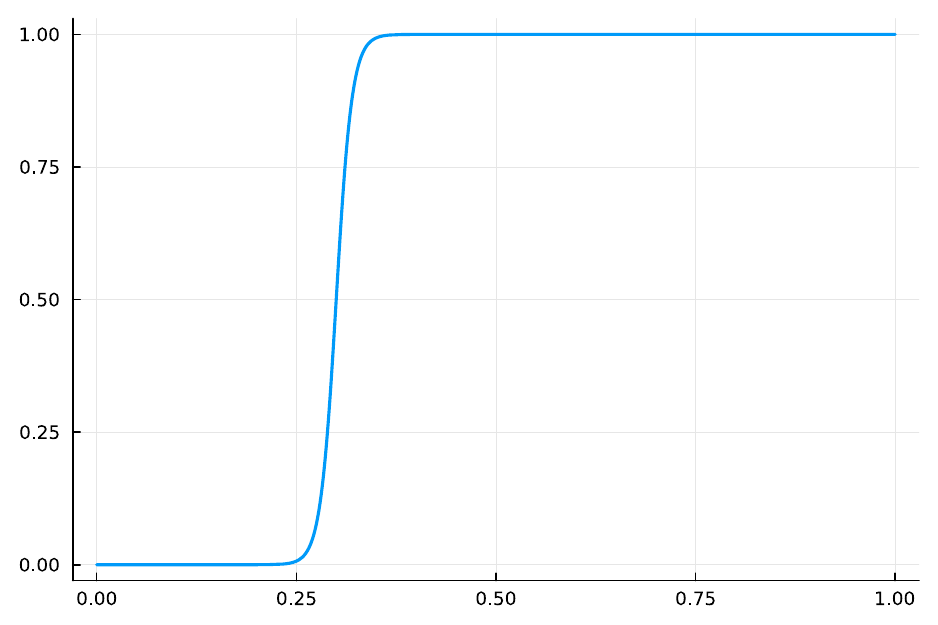}
\includegraphics[height=5cm]{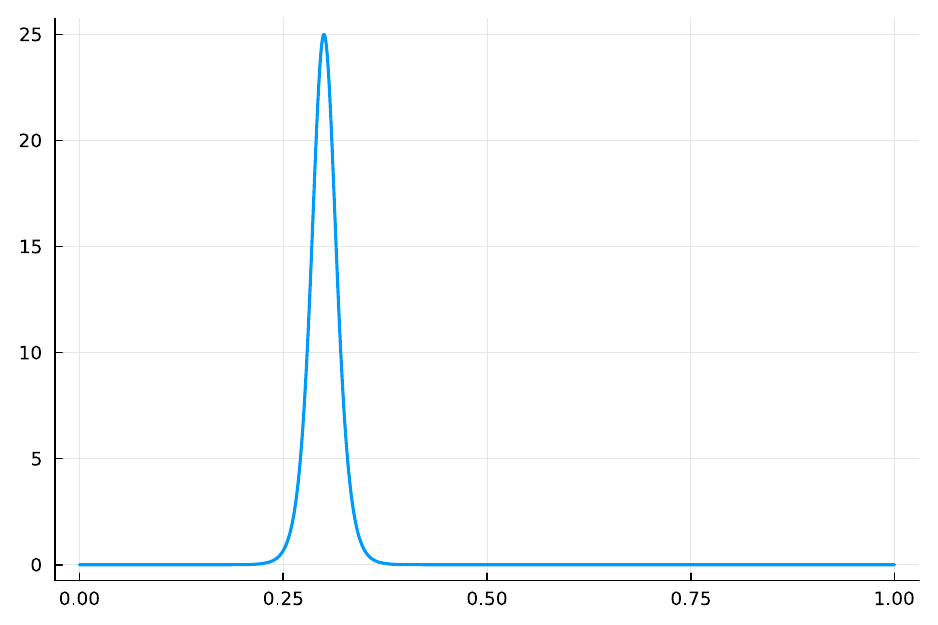}
\caption{The BUM function $Q(t)=\frac{1}{1+\text{e}^{-M(t-\alpha)}}$ for $M=100$, $\alpha=0.3$, and the corresponding weight function $w(t)=\frac{d Q(t)}{d t}$.} \label{figweight}
\end{figure}
Computing $\Sp_Q(\pmb{C}^\top\pmb{x})$ can be hard even if $Q$ is continuously differentiable. Choose $Q(t)=\frac{1}{1+\text{e}^{-M(t-\alpha)}}$ for some large constant $M>0$. We get $Q(0)\simeq 0$ and $Q(1)\simeq 1$. Notice that $Q(t)$ is a continuous approximation of the step function. The corresponding weight function is $w(t)=\frac{M\text{e}^{-M(t-\alpha)}}{(1+\text{e}^{-M(t-\alpha))^2}}$ (see Figure~\ref{figweight}). Then $\Sp_w(\pmb{C}^\top\pmb{x})\simeq\VAR_{\alpha}(\pmb{C}^\top\pmb{x})$. 

\subsubsection{Solution algorithms}
\label{sec:4}

To solve the problem~(\ref{defowaprob}) with  independent and uniformly distributed cost coefficients $C_i,$  $i\in[n]$, 
one can use the  sampling algorithm (Algorithm~\ref{alg:sampling}).
In this section, we propose alternative approximation methods that boil down to solving deterministic counterparts of~(\ref{defowaprob}) with a linear or second order cone objective function. In consequence, if the set of feasible solutions~$\X$ is convex, then we get polynomially solvable convex approximations of~(\ref{defowaprob})~\cite{AG03, BV04} that can be solved efficiently using standard solvers. Furthermore, for a concave BUM function $Q$ we get some approximation guarantees. Notice that~(\ref{defowaprob}) reduces to~(\ref{defopt}) if $Q$ is such that $Q(t)=1-Q(1-t)$ for $t\in [0,1]$. By Lemma~\ref{lem:symmetric}, it is then enough to solve~(\ref{defopt}) for the expected cost vector $\hat{\pmb{c}}=\E[\pmb{C}]$.

Let the vector $\pmb{c}^Q$ be defined as~(\ref{cqdef}). For the random cost $C_i$ uniformly distributed in $[\underline{c}_i,\overline{c}_i]$, we get $c^Q_i=\Sp_{Q}(C_i)=\int_0^1 (\overline{c}_i-t(\overline{c}_i-\underline{c}_i)) d Q(t)=\underline{c}_i+(\overline{c}_i-\underline{c}_i)\int_0^1 Q(t) dt$, using the integration by parts.

\begin{lemma}
\label{lemappruni}
If the BUM function $Q$ is concave and  $\pmb{x}^*\in \X\subseteq \R_{+}^n$ is an optimal solution to~(\ref{defopt}) for the cost vector $\pmb{c}^Q$, then
\begin{equation}
\forall \pmb{x}\in \X\quad \Sp_Q(\pmb{C}^\top\pmb{x}^*)\leq 
\beta_L \cdot\Sp_Q(\pmb{C}^\top\pmb{x}),
\label{apspuni}
\end{equation}
where 
\[
\beta_L=2\int_0^1 Q(t)\, d t
\]
and $\frac{1}{2}\leq \int_0^1 Q(t)\, d t\leq 1$.
\end{lemma}
\begin{proof}
Let $\beta=\int_0^1 Q(t)\, d t$. Since the BUM function $Q$ is concave, $Q(t)\geq t$ for $t\in [0,1]$ and
$\frac{1}{2}\leq \beta\leq 1$. Thus
\begin{align*}
\Sp_Q(\pmb{C}^\top\pmb{x}^*)&\leq \sum_{i\in [n]} x_i^*\Sp_Q(C_i)\leq \sum_{i\in [n]} x_i\Sp_Q(C_i)
=\underline{\pmb{c}}^\top\pmb{x}+\beta(\overline{\pmb{c}}^\top\pmb{x}-\underline{\pmb{c}}^\top\pmb{x})=(1-\beta)\underline{\pmb{c}}^\top\pmb{x}+\beta\overline{\pmb{c}}^\top \pmb{x}\\
& \leq \beta(\underline{\pmb{c}}^\top\pmb{x}+\overline{\pmb{c}}^\top \pmb{x})= 2\beta\E[\pmb{C}^\top\pmb{x}]\leq 2\beta \cdot\Sp_Q(\pmb{C}^\top\pmb{x})=\beta_L \cdot\Sp_Q(\pmb{C}^\top\pmb{x}).
\end{align*}
The first inequality follows from Property~\ref{propdist} (positive homogeneity, subadditivity),
the second follows from the optimality of~$\pmb{x}^*$,
and the last two result from the fact that $\frac{1}{2}\leq \beta\leq 1$ and inequality~(\ref{propbound1}).
\end{proof}
Lemma~\ref{lemappruni} leads to the following theorem:
\begin{theorem}
\label{thmappruni}
Under the assumptions of Lemma~\ref{lemappruni}, the problem~(\ref{defowaprob}) is approximable 
within~$\beta_L$, if $\int_0^1 Q(t)\, d t$  and an optimal solution to the problem~(\ref{defopt}) for~$\pmb{c}^Q$
can be computed in  polynomial time.
\end{theorem}

The constant $\beta_L$ for several BUM functions is shown in Table~\ref{tabappr}. Observe that the problem with the Conditional Value at Risk can be approximated within $2-\alpha$, which is better than $2$ for $\alpha \in (0,1]$. Note also that for the BUM function~(\ref{eqp}) the approximation ratio is less than~$\frac{3}{2}$.

We now show approximation methods for~(\ref{defowaprob})
based on deterministic mathematical programming models with   linear or second-order cone objective functions.
We start with making a reasonable assumption, that
is, let $\gamma$ be a constant such that  $\hat{c}_i\geq \gamma  (\overline{c}_i-\underline{c}_i)$ for each $i\in [n]$.
 Notice that $\gamma\in [\frac{1}{2},\infty)$ and the bound is tight ($\gamma=\frac{1}{2}$) if $\underline{c}_i=0$.
 Let
$$\beta_V=\frac{1}{\sqrt{12}}\int_{0}^1\sqrt{t^{-1}} \, d Q(t)$$
and
$$\beta_H=\frac{\sqrt{2}}{2}\int_{0}^1 \sqrt{\ln t^{-1}}\, d Q(t), $$
where $t\in (0,1)$ and 
 formulate the following optimization problem:
\begin{equation}
\label{convappr1}
\begin{array}{ll}
	\min & \displaystyle\hat{\pmb{c}}^\top \pmb{x}+\beta\sqrt{\sum_{i\in [n]} (\overline{c}_i-\underline{c}_i)^2x_i^2}\\
	& \pmb{x}\in \mathcal{X}\subseteq \R^n_{+}
\end{array}
\end{equation}
where $\beta=\min\{\beta_V,\beta_H\}$. The constants $\beta_V$ and $\beta_H$ for various BUM functions are shown in Table~\ref{tabappr}. 
\begin{table}[ht]
\caption{Constant $\beta_L$, $\beta_V$ and $\beta_H$ for various concave BUM functions.} \label{tabappr}
\centering
\begin{tabular}{l|lll}
BUM function $Q(t)$ & $\beta_L$ & $\beta_V$  & $\beta_H$\\ \hline
$\displaystyle\min\{t/\alpha,1\}$, $\alpha\in (0,1]$ & $2-\alpha$ & $\frac{1}{\sqrt{3\alpha}}$ &  $\frac{\sqrt{2}}{2\alpha}
\Gamma(\frac{3}{2}, -\ln\alpha)$ \\
$\displaystyle\frac{p}{p-1}(t-\frac{1}{p}t^p)$, $p> 1$ & $\displaystyle\frac{p+2}{p+1}$ & $\frac{2p}{\sqrt{3}(2p-1)}$ & $\frac{\sqrt{2\pi}}{2} \cdot \frac{p-(1/\sqrt{p})}{p-1}$\\
$\displaystyle t^\phi$, $\phi \in (0,1]$ & $\displaystyle\frac{2}{\phi+1}$ & $\frac{\phi}{\sqrt{12}(\phi-0.5)}$, $\phi \in (0.5,1]$ & $\frac{\sqrt{\pi}}{\sqrt{2\phi}}$
\end{tabular}
\end{table}

\begin{lemma}
\label{lapprsoc}
If the BUM function $Q$ is concave and $\pmb{x}^*$ is an optimal solution to~(\ref{convappr1}), then $\Sp_Q(\pmb{C}^\top\pmb{x}^*)\leq (1+\beta/\gamma) \cdot \Sp_Q(\pmb{C}^\top\pmb{x})$ for all $\pmb{x}\in \X$.
\end{lemma}
\begin{proof}
See Appendix~\ref{dod}.
\end{proof}
\begin{theorem}
\label{apprsoc}
Under the assumptions of Lemma~\ref{lapprsoc}, the problem~(\ref{defowaprob}) is approximable 
within~$1+\beta/\gamma$, if the value of $\beta$  and an optimal solution to the problem~(\ref{convappr1}) 
can be computed in  polynomial time.
\end{theorem}
Since $\gamma\geq \frac{1}{2}$, the largest bound in Theorem~\ref{apprsoc} is $1+2\beta$. This bound occurs, for example, if all lower bounds of the cost intervals are~0 ($\underline{c}_i=0$ for each $i\in [n]$), which is rather unrealistic in practice. The bound $1+2\beta=1+2\min\{\beta_V, \beta_H\}$ is greater than $\beta_L$ for the BUM functions listed in Table~\ref{tabappr}. However, $1+\beta/\gamma$ can be smaller than $\beta_L$ for some practical instances with larger $\gamma$ (notice that $\gamma$ depends on the input data).

We now propose a modification of~(\ref{convappr1}) which may give better results for practical instances.  Let us define the constants:
$$\beta^{(1)}_B=\frac{1}{3}  \int_0^1  \ln t^{-1}\, d Q(t),$$
$$\beta^{(2)}_B=\frac{1}{\sqrt{6}}  \int_0^1  \sqrt{\ln t^{-1}}\, d Q(t).$$
Consider now the following optimization problem:
\begin{equation}
\label{convapprbern}
\begin{array}{lll}
	\min & \displaystyle \hat{\pmb{c}}^\top \pmb{x}+\beta^{(1)}_B t + \beta^{(2)}_B\sqrt{\sum_{i\in [n]} (\overline{c}_i-\underline{c}_i)^2 x_i^2}\\
	&t\geq (\overline{c}_i-\underline{c}_i)x_i & i\in [n]\\
	& \pmb{x}\in \mathcal{X}\subseteq \R^n_{+}
\end{array}
\end{equation}	

\begin{lemma}
\label{lapprsoc1}
If the BUM function $Q$ is concave and $\pmb{x}^*$ is an optimal solution to~(\ref{convapprbern}), then $\Sp_w(\pmb{C}^\top\pmb{x}^*)\leq (1+(\beta^{(1)}_B+\beta^{(2)}_B)/\gamma) \cdot\Sp_w(\pmb{C}^\top\pmb{x})$ for all $\pmb{x}\in \X$.
\end{lemma}
\begin{proof}
See Appendix~\ref{dod}
\end{proof}
\begin{theorem}
\label{apprsoc1}
Under the assumptions of Lemma~\ref{lapprsoc1}, the problem~(\ref{defowaprob}) is approximable 
within~$(1+(\beta^{(1)}_B+\beta^{(2)}_B)/\gamma) $, if $\beta^{(1)}_B$, $\beta^{(2)}_B$  and an optimal solution to the problem~(\ref{convapprbern}) 
can be computed in  polynomial time.
\end{theorem}
Observe that $\beta^{(2)}_B< \beta_H$. However $\beta^{(1)}_B+\beta^{(2)}_B\geq  \beta_H$, so the theoretical approximation ratio of~(\ref{convapprbern}) is not better than~(\ref{convappr1}). But in practice, the solutions obtained by~(\ref{convapprbern}) can be better, as more information is taken into account to solve this problem.

\begin{remark}
If the BUM function $Q$ is not concave, then the objective functions  of~(\ref{defopt}) for $\pmb{c}^Q$, (\ref{convappr1}) and~(\ref{convapprbern}) are upper bounds on $\Sp_Q(\pmb{C}^\top\pmb{x})$. Therefore, we can still use the models to obtain convex (linear or SOC) approximations of~(\ref{defowaprob}), but without performance guarantees.
\end{remark}

\section{Computational tests}
\label{sec4}

\subsection{Setup}

The purpose of this section is to compare the performance of the algorithms that have been introduced in Section~\ref{sec:4} to the performance of the sampling algorithm (Algorithm~\ref{alg:sampling}). Specifically, the following algorithms are tested:
\begin{itemize}
\item The sampling algorithm with $K\in\{10,50,100\}$, where $K$ scenarios are sampled uniformly and an OWA problem with discrete scenario set is solved.

\item Solving the nominal problem \eqref{defopt} with cost vector $\pmb{c}^Q$ (see Theorem~\ref{thmappruni}). We denote this method as M1.

\item Solving the nonlinear program \eqref{convappr1} (see Theorem~\ref{apprsoc}). We denote this method as M2.

\item Solving the nonlinear program \eqref{convapprbern} (see Theorem~\ref{apprsoc1}). We denote this method as M3.
\end{itemize}

Each solution is then evaluated by sampling a large number of scenarios (in this case, 10000) and calculating the discrete OWA value for this set. We measure both solution quality as well as computation time. Models were solved using Gurobi Version 12.0.0 on a server with an AMD EPYC 9474F 48-Core Processor and 264GB RAM. All processes were restricted to one thread and all Gurobi models were given a time limit of 30 seconds.

\subsection{Instance generation}

For the purpose of this comparison, we generate min-knapsack instances of the type
\[ \min \left\{ \sum_{i\in[n]} c_i x_i : \sum_{i\in[n]} w_i x_i \ge B,\ \pmb{x}\in\{0,1\}^n \right\} \]
in the following way. Item weights $w_i$ are chosen uniformly i.i.d. from the set $\{50,\ldots,100\}$. We then set $B=\sum_{i\in[n]} w_i/10$. To generate cost intervals, we set $\underline{c}_i = (w_i-r_i)(1+\underline{\epsilon}_i)$ and $\overline{c}_i = (w_i + r_i)(1+\overline{\epsilon}_i)$
with $r_i$ chosen  uniformly i.i.d. from $\{0,\ldots,49\}$, while $\underline{\epsilon}_i$ and $\overline{\epsilon}_i$ are chosen uniformly i.i.d. from $[-0.1,0.1]$. If $\underline{c}_i > \overline{c}_i - 1$, these values are sampled again to ensure that all intervals have a length of at least one. Note that this generation procedure ensures that item weights and costs are correlated, and interval midpoints are close to $w_i$. A similar idea was applied in~\cite{BZ80} for generating tested instances of 0-1 Knapsack.

As BUM function, we use $Q(t) = \frac{p}{p-1}(t-\frac{1}{p}^p)$ with $p\in\{1.4, 2, 4\}$ (see Figure~\ref{figsampq}). Furthermore, we chose $n\in\{10\cdot 2^k : k=0,\ldots,7\}$. For each combination of $p$ and $n$, we generate 200 instances. All reported values are averaged over these 200 instances.

\subsection{Results}

\begin{figure}[htbp]
\begin{center}
\subfigure[$p=1.4$\label{fig:obj1}]{\includegraphics[width=0.5\textwidth]{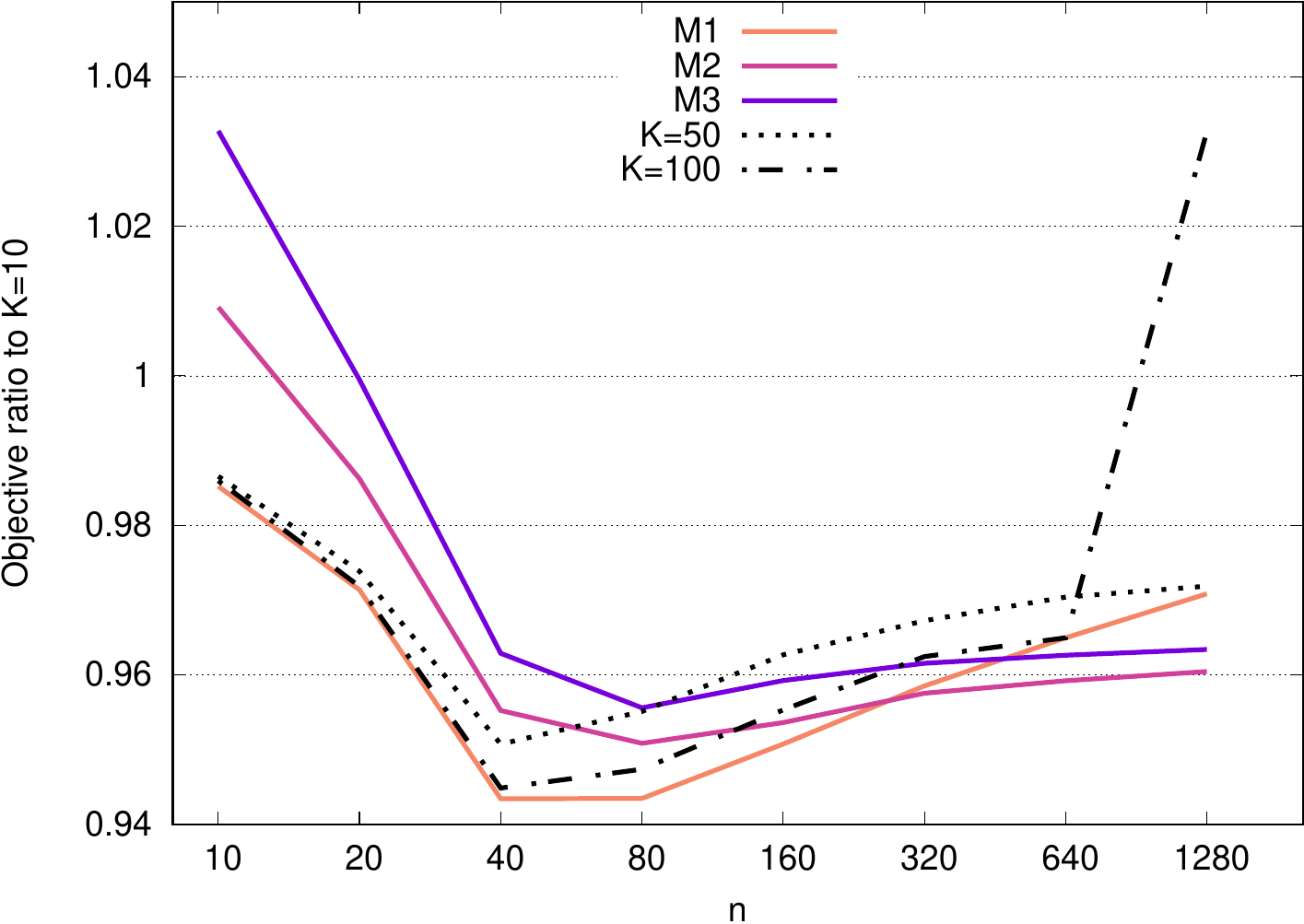}}%
\subfigure[$p=2$\label{fig:obj2}]{\includegraphics[width=0.5\textwidth]{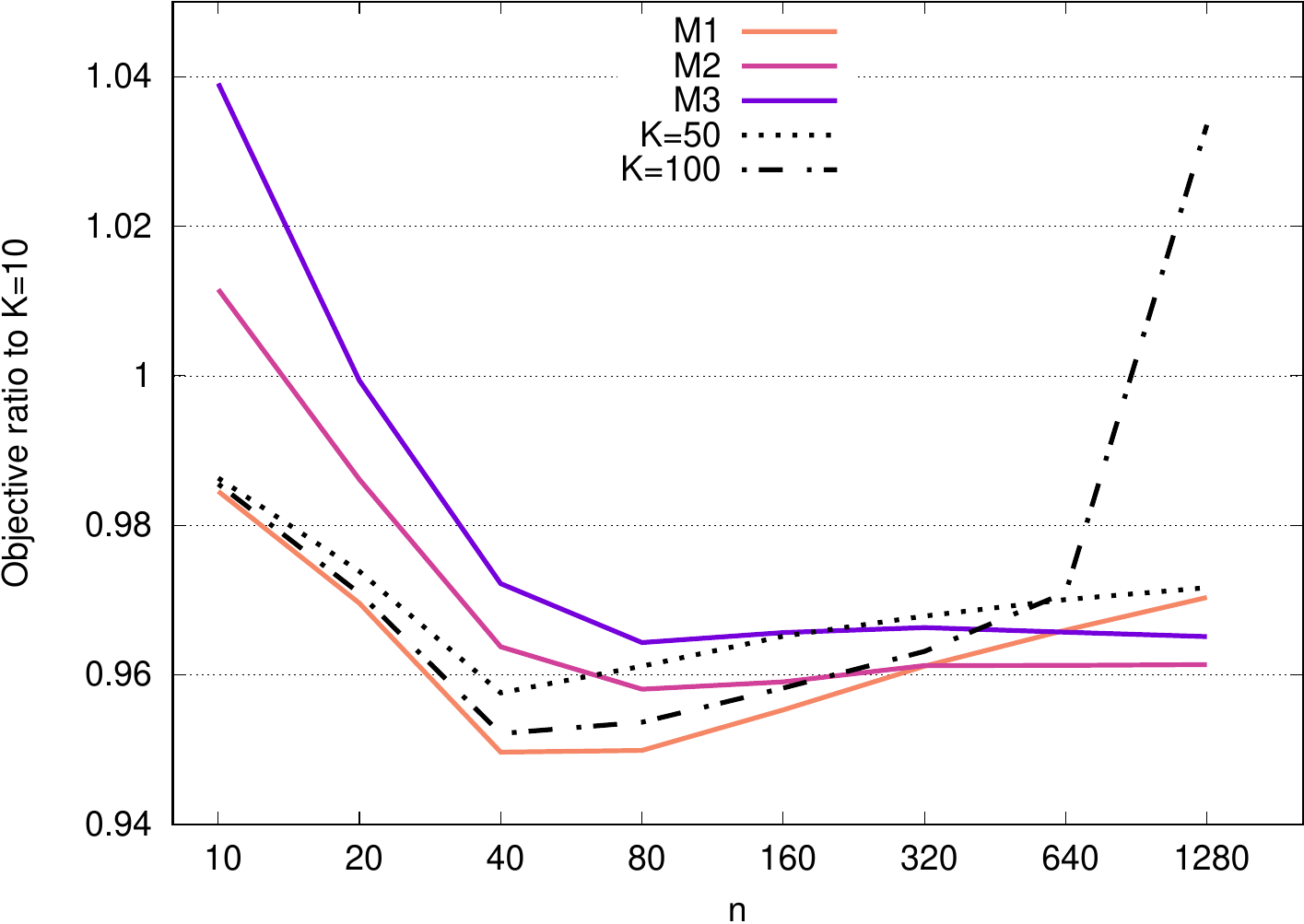}}
\subfigure[$p=4$\label{fig:obj3}]{\includegraphics[width=0.5\textwidth]{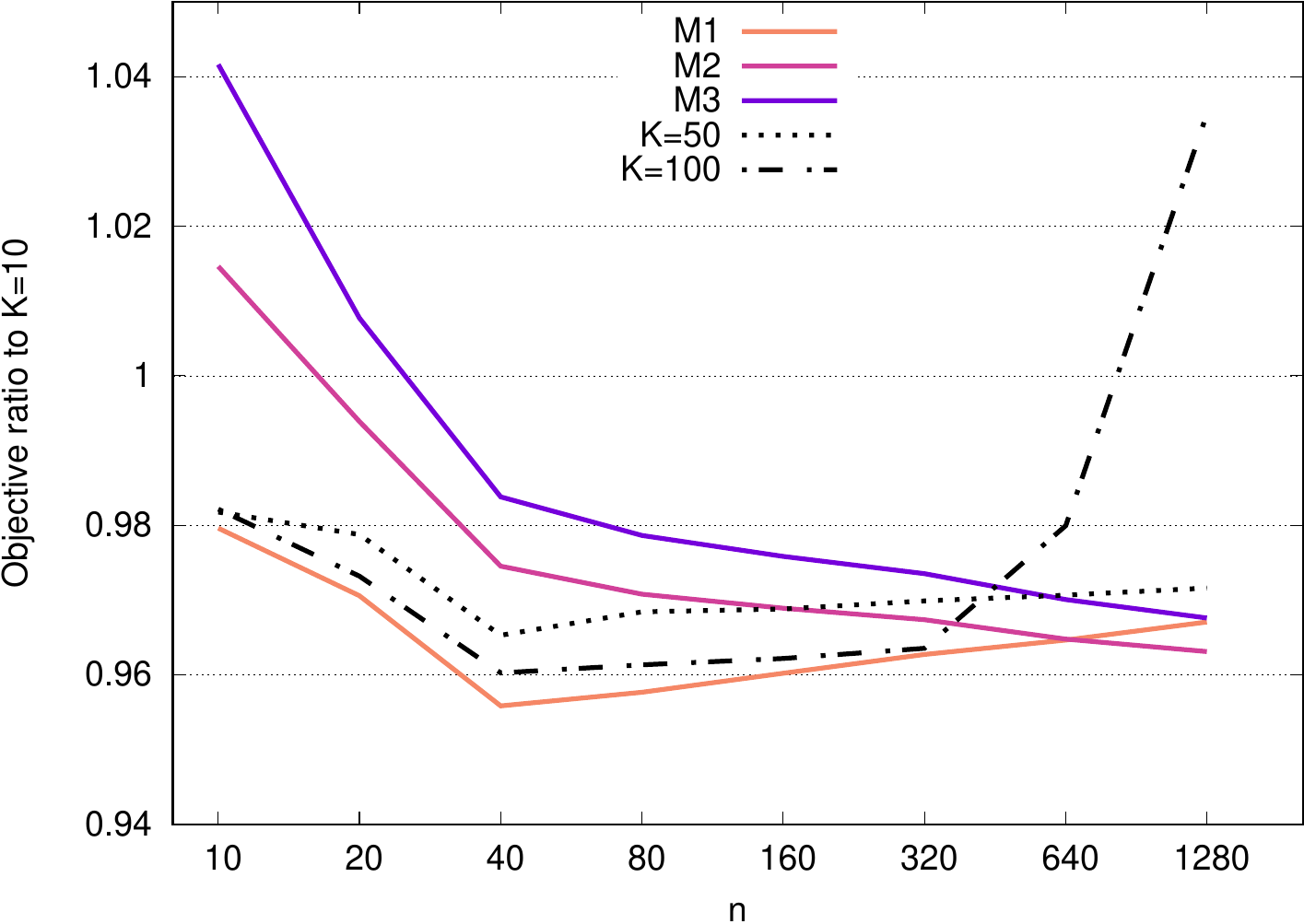}}
\caption{Average ratio of objective value to objective value of sampling with $K=10$.}\label{fig:obj}
\end{center}
\end{figure}

We first assess the quality of solutions that we find, by calculating for each instance the objective value that has been achieved, divided
by the objective value of the sampling method with $K=10$. We therefore use this method as a baseline for our comparison, with values smaller than one indicating better performance than sampling with $K=10$, and values larger than one indicating worse performance. Our averaged results are presented in Figure~\ref{fig:obj}.

First, we note that while results slightly differ between different values of $p$, they remain similar in terms of qualitative insights. We note that with more sampled scenarios ($K=50$ and $K=100$), the estimation of the objective function improves and therefore also the quality of the solutions that we obtain improves. However, this relative advantage declines with increasing $n$. For $n=1280$, the solution with $K=100$ performs even worse than for $K=10$, which is due to the imposed time limit (which we discuss next).

Furthermore, we see that M1 performs well throughout all settings, and consistently outperforms both sampling approaches. As $n$ increases, this advantage becomes smaller; for the largest values of $n$ we studied, M2 and M3 begin to outperform M1. Finally, M2 consistently performs better than M3.

\begin{figure}[htbp]
\begin{center}
\subfigure[$p=1.4$\label{fig:times1}]{\includegraphics[width=0.5\textwidth]{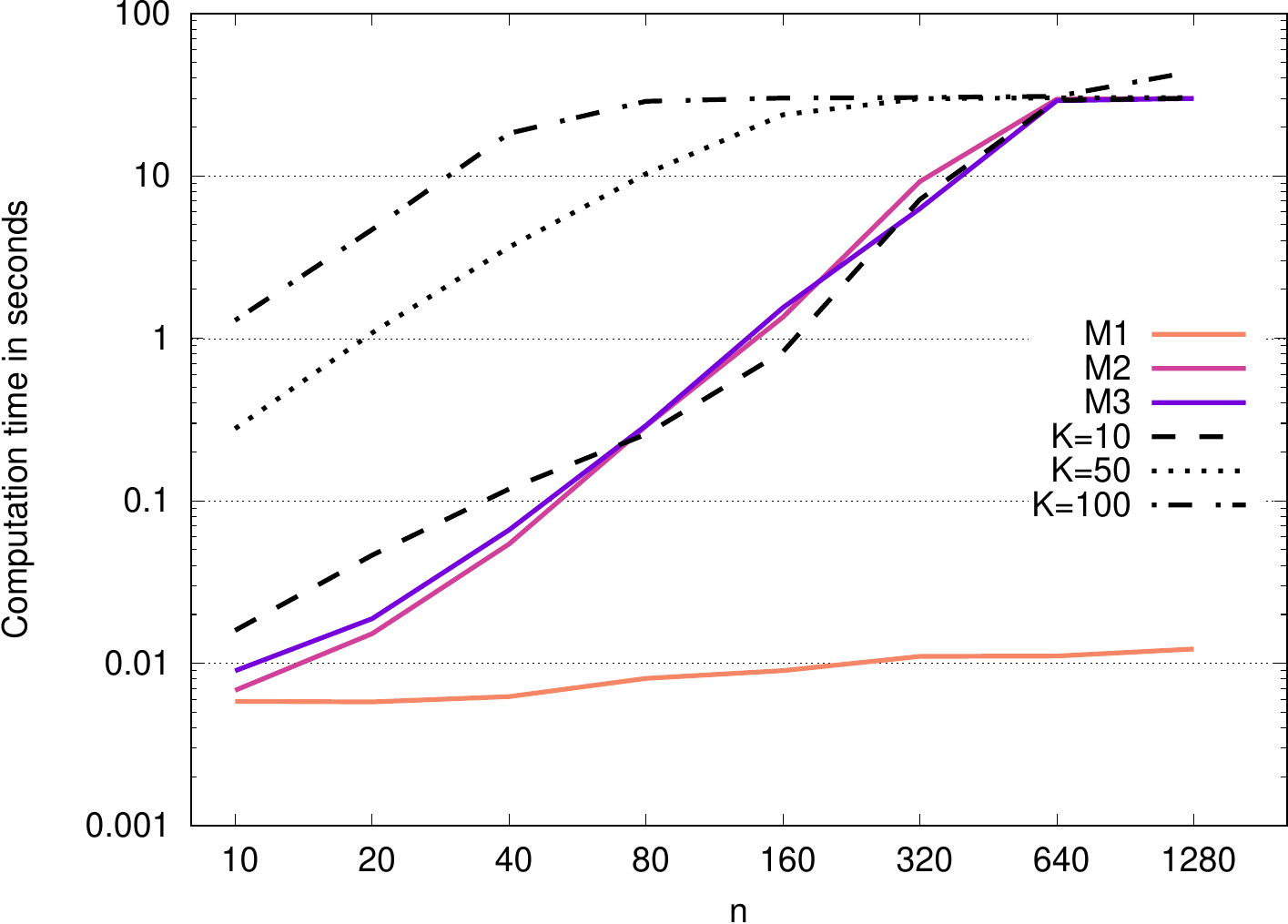}}%
\subfigure[$p=2$\label{fig:times2}]{\includegraphics[width=0.5\textwidth]{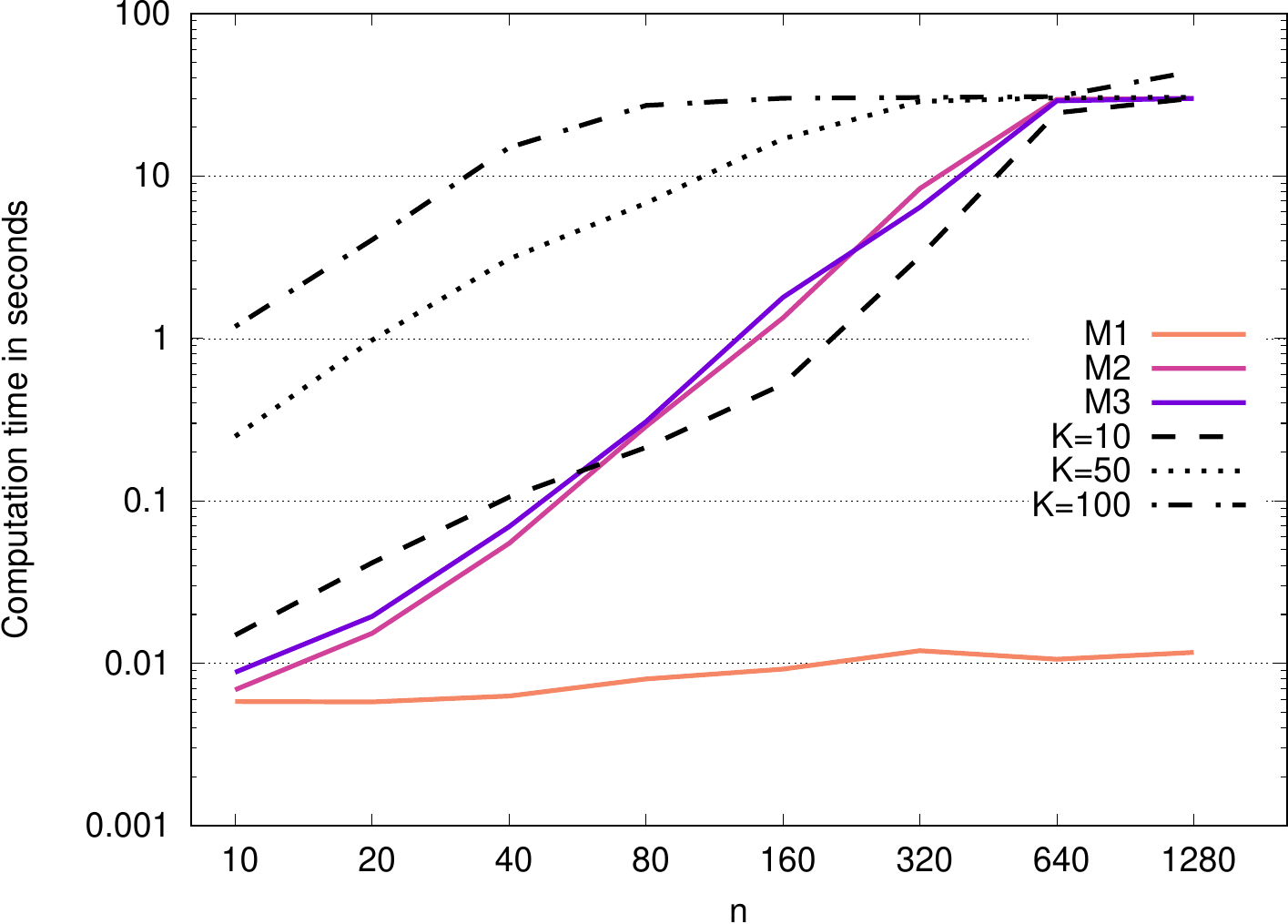}}
\subfigure[$p=4$\label{fig:times3}]{\includegraphics[width=0.5\textwidth]{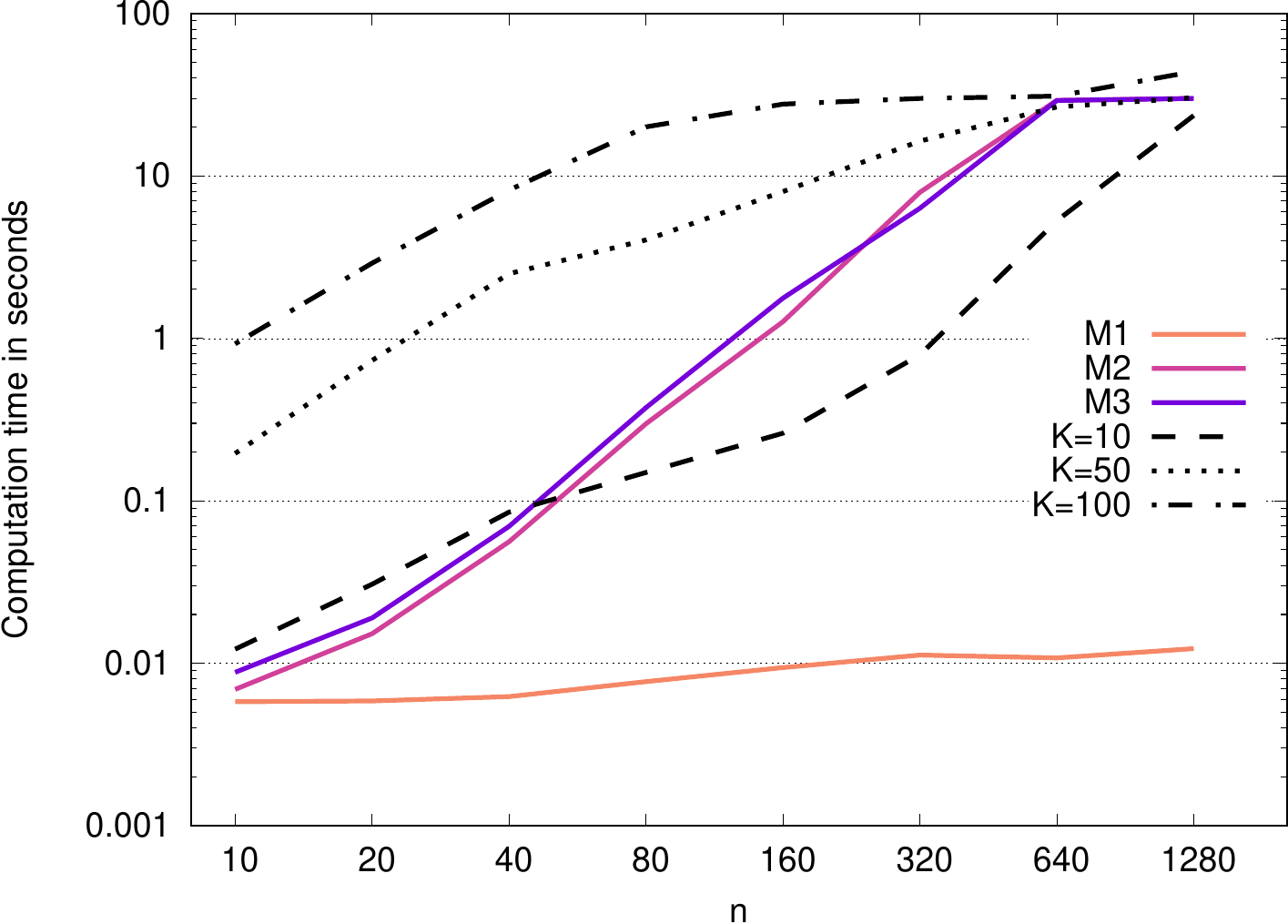}}
\caption{Average computation times.}\label{fig:times}
\end{center}
\end{figure}

We now turn to computation times, which are presented in Figure~\ref{fig:times}. Note the double-logarithmic axes. Computation times for the sampling approach are very sensitive to the number of scenarios $K$. While the baseline $K=10$ remains relatively easy to solve, for $K=100$ the time limit of 30 seconds is limiting the computations already for small values of $n$ (for $n=1280$, Gurobi sometimes exceeded the given time limit, which leads to averages above 40 seconds). Methods M2 and M3 have similar computational effort, which is comparable to the effort of $K=10$, in particular for $p=1.4$.  Method M1, which simple requires a min-knapsack problem to be solved, scales best amongst all methods.

To summarize these findings, we see that the sampling approach can be used to approximate the interval-based problem faithfully, but requires high computational effort. Our proposed alternative methods can achieve a better trade-off between computational effort and solution quality; in particular the simple model M1 performs best in both solution quality and computation times for $n$ up to 320.

\section{Conclusions}

In this paper, we have extended the traditional OWA optimization to the continuous case. For this purpose, we have adopted the concept of a distortion risk measure, whose idea is similar to OWA. Special cases of distortion risk measure are Value at Risk and Conditional Value at Risk, commonly used in stochastic optimization. We have investigated more deeply the model in which the objective function coefficients are independent random variables uniformly distributed in their interval supports. Using an appropriate BUM function, we can model various risk aversions of decision-makers.  Unfortunately, evaluating a given solution turned out to be $\#$P-hard. Therefore, computing an exact solution to the resulting problem can be a difficult task. Notice, however, that the $\#$P-hardness has been shown for non-concave BUM function (modeling VaR) and the computational complexity remains open when BUM is concave. We have shown that the discrete OWA is an estimator of the distortion risk measure. This fact allows us to evaluate approximately a given solution  for a large sample of random cost realizations.  We have proposed some approximation algorithms for solving the problems. For the concave BUM function, we have provided  approximation guarantees for these algorithms. The computational tests performed for 0-1 Knapsack instances suggest that a simple linear approximation, i.e. solving 
the nominal problem \eqref{defopt} with cost vector $\pmb{c}^Q$, 
can perform well in practice.

\subsubsection*{Acknowledgements}
Marc Goerigk and Werner Baak are supported by the Deutsche Forschungsgemeinschaft (DFG) through grant 448792059. Adam Kasperski and Pawe{\l} Zieli{\'n}ski are supported by the National Science Centre, Poland, grant 2022/45/B/HS4/00355.

\appendix

\section{Appendix}
\label{dod}

\begin{proof}[The proof of inequality~(\ref{propbound2})]
We show that $\Sp_{Q}(X) -\E[X]\geq 0$.
The following computation follows from the properties of Riemann-Stieltjes integral~(see, e.g.,~\cite{R76}):
\begin{align*}
\int_0^1 \VAR_{1-t}(X)\, d Q(t)-\E[X]&= \int_0^1 \VAR_{u}(X)\, d (1-Q(1-u))- \int_0^1 \VAR_{u}(X)\, d u\\
&=\int_0^1 \VAR_{u}(X)\, d (1-Q(1-u)-u),
\end{align*}
where the first equality follows changing variable~$t$, $t=1-u$, and the second one from the combination property.
The above integrals are well-defined because $\VAR_{1-t}(X)$ is continuous and $Q(t)$ has bounded variation.
Let $B(u)=1-u -Q(1-u)$, $B(u)$ is  such $B(0)=0$ and $B(1)=1$. Moreover $B(u)\leq 0$ for every $u\in[0,1]$,
since $Q(1-u)\geq 1-u$, which due to the fact that~$Q$ is a concave function.
By the integration by parts formula, we have
\begin{align*}
\int_0^1 \VAR_{u}(X)\, d B(u)&= \VAR_{1}(X)B(1)-\VAR_{0}(X)B(0)
-\int_0^1 B(u)   d\, \VAR_{u}(X)\\
&=\int_0^1(-B(u))   d\, \VAR_{u}(X).
\end{align*}
Clearly, $\int_0^1(-B(u))   d\, \VAR_{u}(X)\geq 0$,  since $-B(u)\geq 0$ and $\VAR_{u}(X)$ is a nondecreasing function.
Thus $\Sp_{Q}(X)\geq \E[X]$.
\end{proof}

\begin{proof}[The proof of inequality~(\ref{propbound2a})]
We choose $\beta\geq 1$ such that $\beta\,\E[X]-\Sp_{Q}(X)\geq 0$. A simple computation gives
\begin{align*}
\beta \int_0^1 \VAR_{1-t}(X)\, d t-\int_0^1 \VAR_{1-t}(X)\, d Q(t)= \int_0^1 \VAR_{1-t}(X)\, d (\beta t-Q(t)).
\end{align*}
Let $B(u)= \beta t-Q(t)$, $B(t)$ is the function such $B(0)=0$ and $B(1)=\beta-1$.
By the integration by parts formula, we have
\begin{align*}
\int_0^1 \VAR_{1-t}(X)\, d B(t)&= \VAR_{0}(X)B(1)-\VAR_{1}(X)B(0)
-\int_0^1 B(u)   d\, \VAR_{1-t}(X)\\
&=\underline{x}(\beta-1)+ \int_0^1B(t)   d\, (-\VAR_{1-t}(X)).
\end{align*}
By  assumption $\underline{x}\geq 0$ and $Q$ is a concave function, i.e. $Q(t)\geq t$ for every $t\in [0,1]$.
Hence and the fact that $-\VAR_{1-t}(X)$ is nondecreasing function,
we conclude that  $\underline{x}(\beta-1)+ \int_0^1B(t)   d\, (-\VAR_{1-t}(X))\geq 0$, if
$\beta\geq\sup_{t\in (0,1]}\frac{Q(t)}{t}$.
Combining this with  the inequality~(\ref{propbound2}) gives  the inequalities~(\ref{propbound2a}).
\end{proof}

\begin{proof}[The proof of equality~(\ref{esym})]
Let $\mu=\E[X]$. Since $X$ has a symmetric probability distribution,
$1-F_{X}(y+\mu)=F_{X}(\mu-y)$ for every $y\in \R$. Hence and
by the continuity of $\VAR_{t}$ for $t\in [0,1]$, we have
$\VAR_{1-t}(X)=y+\mu \Leftrightarrow \VAR_{t}(X)=-y+\mu$.
Thus,
\begin{equation}
\label{f0}
\VAR_{1-t}(X)+ \VAR_{t}(X)=2\mu.
\end{equation}
By assumption $Q$ is such that $Q(t)=1-Q(1-t)$ for $t\in[0,1]$. 
From this and~(\ref{f0}), we get
\begin{align*}
\Sp_{Q}(X)&=\int_0^1 \VAR_{1-t}(X)\, d Q(t)=\int_0^1 (2\mu -\VAR_{t}(X))\, d Q(t)\\
&=
2\mu-\int_0^1\VAR_{t}(X)\, d(1-Q(1-t))=2\mu-\Sp_{Q}(X).
\end{align*}
Therefore, $\Sp_{Q}(X)=\mu$.
\end{proof}

The next proofs are based on the following well-known inequalities
 (see, e.g.,~\cite{BLB03}).
The following bounds on~the value at risk  results from Chebyshev’s inequality.
\begin{lemma}
\label{chbound}
Let  $X$  be a random variable  with bounded support. Then  for any  $t\in (0,1]$,
\begin{equation}
\label{chappr}
\VAR_{1-t}(X)\leq  \E[X]+\frac{1}{\sqrt{t}}\sqrt{\V[X]}.
\end{equation}
\end{lemma}

\begin{lemma}[Hoeffding's inequality]
Let $X_1,\ldots,X_n$ be independent random variables such that $a_i\leq X_i\leq b_i$, $ i\in [n]$,
and $S_n=\sum_{i\in [n]} X_i$.
Then for any $\delta>0$,
\[
\Pr[ S_n- \E[S_n] \geq \delta]\leq \mathrm{exp}\left(-\frac{2\delta^2}{\sum_{i\in [n]}(b_i-a_i)^2}  \right).
\]
\label{hoeinq}
\end{lemma}
\begin{lemma}[Bernstein's inequality]
Let $X_1,\ldots,X_n$ be independent zero mean random variables such that $|X_i|\leq M_i$, $ i\in [n]$,
and $S_n=\sum_{i\in [n]} X_i$.
Then for any $\delta>0$ the following inequality
\[
\Pr[ S_n \geq \delta]\leq \mathrm{exp}\left(-\frac{\delta^2}{2\sum_{i\in [n]}\V[X_i]+\frac{2}{3} M_{\max}\delta}  \right) 
\]
holds, where $\V[X_i] $ is the variance of~$X_i$ and $M_{\max}= \max_{i\in [n]} M_i$.
\label{berninq}
\end{lemma}

\begin{proof}[The proof of Lemma~\ref{lapprsoc}]
We have $\E[\pmb{C}^\top\pmb{x}]=\hat{\pmb{c}}^\top \pmb{x}$. Because the components of $\pmb{C}$  are independent, we also have $\V[\pmb{C}^\top\pmb{x}]=\sum_{i\in [n]} \V[C_i x_i]=\sum_{i\in [n]} \V[C_i]x_i^2=\frac{1}{12}\sum_{i\in [n]} (\overline{c}_i-\underline{c}_i)^2x_i^2$. Using Lemma~\ref{chbound}, we get the following upper bound on $\Sp_Q(\pmb{C}^\top\pmb{x})$:
$$\Sp_{Q}(\pmb{C}^\top\pmb{x})\leq\int_0^1 \left(\hat{\pmb{c}}^\top \pmb{x}+\frac{1}{\sqrt{12t}}\sqrt{\sum_{i\in [n]} (\overline{c}_i-\underline{c}_i)^2x_i^2}\right) \, dQ(t) = \hat{\pmb{c}}^\top \pmb{x}+\beta_V\sqrt{\sum_{i\in [n]} (\overline{c}_i^2-\underline{c}_i^2)x_i^2}.$$
Therefore, for every $\pmb{x}\in \X$ we have:
$$
\begin{array}{ll}
\Sp_Q(\pmb{C}^\top \pmb{x}^*) & \leq \displaystyle\hat{\pmb{c}}^\top \pmb{x}^*+\beta_V\sqrt{\sum_{i\in [n]} (\overline{c}_i-\underline{c}_i)^2 x_i^{*2}}\leq \hat{\pmb{c}}^\top \pmb{x}+\beta_V\sqrt{\sum_{i\in [n]} (\overline{c}_i-\underline{c}_i)^2x_i^2}\\
& \displaystyle\leq \hat{\pmb{c}}^\top \pmb{x}+\beta_V\sum_{i\in [n]} \sqrt{(\overline{c}_i-\underline{c}_i)^2}x_i\leq \hat{\pmb{c}}^\top \pmb{x}+(\beta_V/\gamma) \hat{\pmb{c}}^\top \pmb{x}\\
&=(1+\beta_V/\gamma)\hat{\pmb{c}}^\top\pmb{x}\leq (1+\beta_V/\gamma)\Sp_Q(\pmb{C}^\top\pmb{x}).
\end{array}$$
The random variable $\pmb{C}^\top\pmb{x}$ is the sum of independent random variables~$C_ix_i$, $\in [n]$,
where $C_i \sim  U[\underline{c}_i, \overline{c}_i]$. Hence $\underline{c}_i x_i\leq C_ix_i\leq \overline{c}_i x_i$, $i\in [n]$.
We give an upper bound on $\VAR_{1-t}(\pmb{C}^\top\pmb{x})$ using the Hoeffding's inequality. Thus, we need to
find~$y^*$ such that
\[
\Pr[\pmb{C}^\top\pmb{x}\leq y^*]\geq 1-t, \; t\in (0,1),
\]
which equivalent to $\Pr[\pmb{C}^\top\pmb{x}> y^*]\leq t$.
Let $y^*=\E[\pmb{C}^\top\pmb{x}]+\delta^*$.
We now determine~$\delta^*$. By Lemma~\ref{hoeinq}, we have
\[
\Pr[\pmb{C}^\top\pmb{x}>\E[\pmb{C}^\top\pmb{x}]+\delta^*]=
\Pr[\pmb{C}^\top\pmb{x}-\E[\pmb{c}^\top\pmb{x}]\geq\delta^*]\leq 
\mathrm{exp}\left(-\frac{2(\delta^*)^2}{\sum_{i\in [n]}(\overline{c}_i-\underline{c}_i)^2 x^2_i}  \right). 
\]
We choose $\delta^*>0$ such that 
\[
\mathrm{exp}\left(-\frac{2(\delta^*)^2}{\sum_{i\in [n]}(\overline{c}_i-\underline{c}_i)^2 x^2_i}  \right)\leq t
\]
holds. Thus
\[
\frac{2(\delta^*)^2}{\sum_{i\in [n]}(\overline{c}_i-\underline{c}_i)^2 x^2_i}\geq \ln t^{-1}.
\]
Finally, we obtain
\[
\delta^*=\frac{\sqrt{2}}{2}\sqrt{\ln t^{-1}} \sqrt{\sum_{i\in [n]} (\overline{c}_i-\underline{c}_i)^2x_i^2},
\]
which yields 
\[
y^*=\hat{\pmb{c}}^\top\pmb{x}+\frac{\sqrt{2}}{2}\sqrt{\ln t^{-1}} \sqrt{\sum_{i\in [n]} (\overline{c}_i-\underline{c}_i)^2x_i^2}.
\]
Hence, we get an upper bound
\[
\VAR_{1-t}(\pmb{C}^\top\pmb{x})\leq y^*= \hat{\pmb{c}}^\top\pmb{x}+\frac{\sqrt{2}}{2}\sqrt{\ln t^{-1}} \sqrt{\sum_{i\in [n]} (\overline{c}_i-\underline{c}_i)^2x_i^2}.
\]
and
$$\Sp_{Q}(\pmb{x})\leq \hat{\pmb{c}}^\top \pmb{x}+\beta_H\sqrt{\sum_{i\in [n]} (\overline{c}_i^2-\underline{c}_i^2)x_i^2}.$$
The rest of the proof is the same as for  the constant $\beta_V$.
\end{proof}

\begin{proof}[The proof of Lemma~\ref{lapprsoc1}]
Consider a random variable $\pmb{C}^\top\pmb{x}$, which is the sum of independent random variables~$C_ix_i$, $\in [n]$,
where $C_i \sim  U[\underline{c}_i, \overline{c}_i]$.
We  first give an upper bound on $\VAR_{1-t}(\pmb{C}^\top\pmb{x})$. Thus, we need to
find~$y^*$ such that
\[
\Pr[\pmb{C}^\top\pmb{x}\leq y^*]\geq 1-t, \; t\in (0,1),
\]
which equivalent to $\Pr[\pmb{C}^\top\pmb{x}> y^*]\leq t$.
Let $y^*=\E[\pmb{C}^\top\pmb{x}]+\delta^*$ and so
\[
\Pr[\pmb{C}^\top\pmb{x}>\E(\pmb{C}^\top\pmb{x})+\delta^*]=
\Pr[\pmb{C}^\top\pmb{x}-\E(\pmb{C}^\top\pmb{x})\geq\delta^*].
\]
Consider the random variable $\pmb{C}^\top\pmb{x}-\E[\pmb{C}^\top\pmb{x}]$ and denote it by $S_n$.
Clearly $S_n$ is the sum of $n$ independent random variables $X_i$, $i\in[n]$, such that
$X_i \in [\underline{c}_i x_i-\E[C_i]x_i, \overline{c}_i x_i-\E[C_i]x_i]$, $\E[X_i]=0$,
$\V[X_i]=\frac{(\overline{c}_i-\underline{c}_i)^2x_i^2}{12}$ and 
$M_i=\overline{c}_i x_i-\E[C_i]x_i=\frac{(\overline{c}_i-\underline{c}_i)x_i}{2}$.
 By Lemma~\ref{berninq}, we have
\[
\Pr[\pmb{C}^\top\pmb{x}-\E[\pmb{C}^\top\pmb{x}]\geq\delta^*]\leq 
\mathrm{exp}\left(-\frac{(\delta^*)^2}{\frac{\sum_{i\in [n]}(\overline{c}_i-\underline{c}_i)^2 x^2_i}{6} +\frac{2}{3}\delta^* M_{\max}}  \right).
\]
We now determine~$\delta^*>0$ such that
\[
\mathrm{exp}\left(-\frac{(\delta^*)^2}{\frac{\sum_{i\in [n]}(\overline{c}_i-\underline{c}_i)^2 x^2_i}{6} +\frac{2}{3}\delta^* M_{\max}}  \right)
\leq t.
\]
An easy computation shows that
\[
(\delta^*)^2-B(t)\delta^*-C(t)\geq 0,
\]
where $B(t)=\frac{2}{3}M_{\max}\ln t^{-1}$ and $C(t)=\frac{\sum_{i\in [n]}(\overline{c}_i-\underline{c}_i)^2 x^2_i}{6}\ln t^{-1}$
and $M_{\max}=\frac{1}{2}\max_{i\in[n]}\{(\overline{c}_i-\underline{c}_i)x_i\}$.
The equation $(\delta^*)^2-B(t)\delta^*-C(t)= 0$ has the following positive solution
\[
\delta^*=\frac{B(t)+\sqrt{(B(t))^2+4C(t)}}{2}.
\]
Therefore, we obtain an upper bound
\[
\VAR_{1-t}(\pmb{C}^\top\pmb{x})\leq y^*= \E[\pmb{C}^\top\pmb{x}]+\frac{B(t)+\sqrt{(B(t))^2+4C(t)}}{2}\leq
\E[\pmb{C}^\top\pmb{x}]+B(t)+\sqrt{C(t)}
\]
and an upper bound on $\Sp_Q(\pmb{C}^\top\pmb{x})$, namely
\begin{align*}
\Sp_Q(\pmb{C}^\top \pmb{x})&=\int_0^1\VAR_{1-t}(\pmb{c}^\top\pmb{x})\, \mathrm{d} Q(t)
\leq \E[\pmb{C}^\top\pmb{x}]+\int_0^1B(t)\, \mathrm{d} Q(t)+ \int_0^1 \sqrt{C(t)}\, \mathrm{d} Q(t)\\
& =\hat{\pmb{c}}^\top \pmb{x}+\beta^{(1)}_B\max_{i\in[n]}\{(\overline{c}_i-\underline{c}_i)x_i\}+ \beta^{(2)}_B\sqrt{\sum_{i\in [n]} (\overline{c}_i-\underline{c}_i)^2 x_i^2}.
\end{align*}
We get for each $\pmb{x}\in \X$
\begin{align*}
\Sp_Q(\pmb{C}^\top \pmb{x}^*)&\leq \hat{\pmb{c}}^\top \pmb{x}^*+\beta^{(1)}_B\max_{i\in[n]}\{(\overline{c}_i-\underline{c}_i)x^*_i\}+ \beta^{(2)}_B\sqrt{\sum_{i\in [n]} (\overline{c}_i-\underline{c}_i)^2 x_i^{*2}} \\
&\leq \hat{\pmb{c}}^\top \pmb{x}+\beta^{(1)}_B\max_{i\in[n]}\{(\overline{c}_i-\underline{c}_i)x_i\}+ \beta^{(2)}_B\sqrt{\sum_{i\in [n]} (\overline{c}_i-\underline{c}_i)^2 x_i^{2}}\\
&\leq \hat{\pmb{c}}^\top \pmb{x}+\beta^{(1)}_B\max_{i\in[n]}\{(\overline{c}_i-\underline{c}_i)x_i\}+ \beta^{(2)}_B\sum_{i\in [n]} (\overline{c}_i-\underline{c}_i) x_i\\
&\leq \hat{\pmb{c}}^\top \pmb{x}+ (\beta^{(1)}_B+\beta^{(2)}_B)\sum_{i\in [n]} (\overline{c}_i-\underline{c}_i) x_i \\
&\leq \hat{\pmb{c}}^\top\pmb{x}+(\beta^{(1)}_B+\beta^{(2)}_B)/\gamma \hat{\pmb{c}}^\top\pmb{x}\leq (1+(\beta^{(1)}_B+\beta^{(2)}_B)/\gamma) \Sp_Q(\pmb{C}^\top\pmb{x}).
\end{align*}
\end{proof}

\end{document}